\documentclass{amsart}
\usepackage{graphicx}
\usepackage{amssymb}
\newtheorem{thmm}{Theorem}
\newtheorem{thm}{Theorem}[section]

\newtheorem{lem}[thm]{Lemma}
\newtheorem{prop}[thm]{Proposition}
\theoremstyle{definition}
\newtheorem{defn}[thm]{Definition}
\newtheorem{rem}[thm]{Remark}
\newtheorem*{problem}{Frobenius Coin Problem}
\newtheorem*{defn*}{Definition}
\newtheorem*{rems*}{Remarks}
\newtheorem*{rem*}{Remark}

\newtheorem*{nt}{Notation}

\numberwithin{equation}{section}

\def \local-algebra {\Lambda ^0(\mathbb R^2)/(\nabla H)}

\begin{document}

\title[Local symplectic algebra] {Local symplectic algebra of quasi-homogeneous curves}
\author{Wojciech Domitrz}
\address{Warsaw University of Technology\\
Faculty of Mathematics and Information Science\\
Plac Politechniki 1\\
00-661 Warsaw\\
Poland\\
and
Institute of Mathematics\\
Polish Academy of Sciences\\
Sniadeckich 8\\
P.O. Box 137\\
00-950 Warsaw\\
Poland }

\email{domitrz@mini.pw.edu.pl}

\thanks{The work of the author was supported by Institute of Mathematics, Polish Academy of Sciences.
}



\begin{abstract}
We study the local symplectic algebra of parameterized curves
introduced by V. I. Arnold in \cite{Ar1}. We use the method of
algebraic restrictions  to classify symplectic singularities of
quasi-homogeneous curves. We prove that the space of algebraic
restrictions of closed $2$-forms to the germ of a $\mathbb
K$-analytic curve is a finite dimensional vector space. We also
show that the action of local diffeomorphisms preserving the
quasi-homogeneous curve on this vector space is determined by the
infinitesimal action of liftable vector fields. We apply these
results to obtain the complete symplectic classification of curves
with the semigroups $(3,4,5)$, $(3,5,7)$, $(3,7,8)$.
\end{abstract}

\maketitle
\section{Introduction} We study the problem of classification of
parameterized curve-germs in a symplectic space $(\mathbb
K^{2n},\omega)$ up to the symplectic equivalence (for $\mathbb
K=\mathbb R$ or $\mathbb C$). The {\bf symplectic equivalence} is
a right-left equivalence (or $\mathcal A$-equivalence) in which
the left diffeomorphism-germ is a symplectomorphism of $(\mathbb
K^{2n},\omega)$ i. e. it preserves the given symplectic form
$\omega$ in $\mathbb K^{2n}$.

The problem of $\mathcal A$-classification of singularities of
parameterized curves-germs was studied by J. W. Bruce and T. J.
Gaffney, C. G. Gibbson and C. A. Hobbs.  Bruce and Gaffney
(\cite{BG}) classified the $\mathcal A$-simple plane curves and in
\cite{GH} the classification of the $\mathcal A$-simple space
curves was given. The singularity (an $\mathcal A$-equivalence
class) is called simple if it has a neighbourhood intersecting
only finite number of singularities. V. I. Arnold (\cite{Ar2})
classified stably simple singularities of curves. The singularity
is stably simple if it is simple and remains simple after
embedding into a larger space.

 The main tool and the invariant separating the
singularities in $\mathcal A$-classification of curves is the
semigroup of a curve singularity $t\mapsto
f(t)=(f_1(t),\cdots,f_m(t))$(see \cite{GH} and \cite{Ar2}). It is
the subsemigroup of the additive semigroup of natural numbers
formed by the orders of zero at the origin of all linear
combinations of the products of $f_i(t)$.

In \cite{Ar1} V. I. Arnold discovered new symplectic invariants of
parameterized curves. He proved that the $A_{2k}$ singularity of a
planar curve (the orbit with respect to standard $\mathcal
A$-equivalence of parameterized curves) split into exactly $2k+1$
symplectic singularities (orbits with respect symplectic
equivalence of parameterized curves). Arnold posed a problem of
expressing  these invariants in terms of the local algebra's
interaction with the symplectic structure. He proposed to call
this interaction {\bf local symplectic algebra}.

In \cite{IJ1} G. Ishikawa and S. Janeczko classified symplectic
singularities of curves in the $2$-dimensional symplectic space.
All simple curves in this classification are quasi-homogeneous.

Symplectic singularity is {\bf stably simple} if it is simple and
remains simple if the ambient symplectic space is symplectically
embedded (i.e. as a symplectic submanifold) into a larger
symplectic space. In \cite{K} P. A. Kolgushkin classified the
stably simple symplectic singularities of curves (in the $\mathbb
C$-analytic category). All stably simple symplectic singularities
of curves are quasi-homogeneous too.

 In \cite{DJZ2} new symplectic invariants of singular
quasi-homogeneous  subsets of a symplectic space were explained by
the algebraic restrictions of the symplectic form to these
subsets.

The algebraic restriction is a equivalence class of the following
relation on the space of differential $k$-forms:

Differential $k$-forms $\omega_1$ and $\omega_2$ have the same
{\bf algebraic restriction} to a subset $N$ if
$\omega_1-\omega_2=\alpha+d\beta$, where $\alpha$ is a $k$-form
vanishing on $N$ and $\beta$ is a $(k-1)$-form vanishing on $N$.

The algebraic restriction of a $k$-form $\omega_1$ to a subset
$N_1$ and the algebraic restriction of a $k$-form $\omega_2$ to a
subset $N_2$ are {\bf diffeomorphic} if there exists a
diffeomorphism $\Phi$ of $\mathbb K^m$ which maps $N_1$ to $N_2$
such that $\Phi^{\ast} \omega_2$ and $\omega_1$ have the same
algebraic restriction to $N_1$ (for details see section
\ref{def-alg-rest}).

 The results in \cite{DJZ2} were obtained by the following
generalization of Darboux-Givental theorem.

\begin{thmm}[\cite{DJZ2}]\label{Darboux}
Quasi-homogeneous subsets of a symplectic manifold $(M,\omega)$
are locally symplectomorphic if and only if algebraic restrictions
of the symplectic form $\omega$ to them are locally diffeomorphic.
\end{thmm}

This theorem reduces the problem of symplectic classification of
quasi-homo\-ge\-neous subsets to the problem of classification of
algebraic restrictions of symplectic forms to these subsets.

In \cite{DJZ2} the method of algebraic restrictions is applied to
various classification problems in a symplectic space. In
particular the complete symplectic classification of classical
A-D-E singularities of planar curves is obtained, which contains
Arnold's symplectic classification of $A_{2k}$ singularity.

In this paper we return to Arnold's original problem of local
symplectic algebra of a parameterized curve. We show that the
method of algebraic restrictions is a very powerful classification
tool for quasi-homogeneous parameterized curves. This is due to
the several reasons. The most important one is that the space of
algebraic restrictions of germs of closed $2$-forms to a $\mathbb
K$-analytic parameterized curve is a finite dimensional vector
space. This fact follows from the following  more general result
conjectured in \cite{DJZ2}, which we prove in this paper.

\begin{thmm}
\label{main-alg} Let $C$ be the germ of a $\mathbb K$-analytic
curve. Then the space of algebraic restrictions of germs of closed
$2$-forms to $C$ is a finite dimensional vector space.
\end{thmm}

By a {\bf $\mathbb K$-analytic curve} we understand a subset of
$\mathbb K^m$ which is locally diffeomorphic to a $1$-dimensional
(possibly singular) $\mathbb K$-analytic subvariety of $\mathbb
K^m$. Germs of $\mathbb K$-analytic parameterized curves can be
identified with germs of irreducible $\mathbb K$-analytic curves.

The tangent space to the orbit of an algebraic restriction $a$ to
the germ $f$ of a parameterized curve is given by the Lie
derivative of $a$ with respect to germs of liftable vector fields
over $f$. We say that the germ $X$ of a liftable vector field acts
{\bf trivially} on the space of algebraic restriction if the Lie
derivative of any algebraic restriction with respect $X$ is zero.

\begin{thmm}
\label{main-lift} The space of germs of liftable vector fields
over the germ of a parameterized quasi-homogeneous curve which act
nontrivially on the space of algebraic restrictions of closed
$2$-forms is a finite dimensional vector space.
\end{thmm}

Theorem \ref{main-alg} is proved in section \ref{proof}. In
section \ref{qh-alg} we prove Theorem  \ref{main-lift} using the
quasi-homogeneous grading on the space of algebraic restrictions.
We show that there exist quasi-homogeneous bases of the space of
algebraic restrictions of closed 2-forms and of the space of
liftable vector fields which act nontrivially on the space of
algebraic restrictions to a quasi-homogeneous parameterized curve.
These bases are allowed us to prove Theorem \ref{action} that
states that the linear action on the space of algebraic
restrictions of closed $2$-forms to the germ of a
quasi-homogeneous parameterized curve by Lie derivatives with
respect to liftable vector fields determines the action on this
space by local diffeomorphisms preserving this germ of the curve.

 Both the space of algebraic
restrictions of symplectic forms and this linear action are
determined by the semigroup of the curve singularity.

We apply the method of algebraic restrictions and results of
Section \ref{qh-alg} to obtain the complete symplectic
classification of curves with the semigroups $(3,4,5)$, $(3,5,7)$
and $(3,7,8)$ in Sections \ref{345}, \ref{357} and \ref{378}.  The
classification results are presented in Table \ref{ss345}, Table
\ref{ss357} and Table \ref{ss378}. All normal forms are given in
the canonical coordinates $(p_1,q_1,\cdots, p_n,q_n)$ in the
symplectic space $(\mathbb R^{2n},\sum_{i=1}^n dp_i\wedge dq_i)$.
The parameters $c, c_1, c_2$ are moduli. The different singularity
classes are distinguished by discrete symplectic invariants:  the
symplectic multiplicity $\mu_{sympl}(f)$, the index of isotropness
$i(f)$ and the Lagrangian tangency order $Lt(f)$, which are
considered in Section \ref{discrete}

We consider only quasi-homogeneous parameterized curves in this
paper. But there are $\mathcal A$-simple singularities of curves
which are not quasi-homogeneous. For example the curve
$f(t)=(t^3,t^7+t^8)$ is not quasi-homogeneous. Then Theorem
\ref{Darboux} cannot be applied for such curves. But there exists
a generalization of this theorem to any subsets $N$ of $\mathbb
K^m$ (\cite{DJZ2}, section 2.6). In general  there is one more
invariant for the symplectic classification problem which can be
represented as a cohomology class in the second cohomology group
of the complex of $2$-forms with zero algebraic restrictions to
$N$. This cohomology groups vanish for quasi-homogeneous subsets
(\cite{DJZ1}). They are finite dimensional for $\mathbb
C$-analytic varieties  with an isolated singularity (\cite{BH}).
It implies that they are finite dimensional for non
quasi-homogeneous $\mathbb C$-analytic curves. The space of
algebraic restriction of closed $2$-forms to a $\mathbb
K$-analytic curve is finite dimensional too by Theorem
\ref{main-alg}.  But the description of the action of
diffeomorphisms preserving a non quasi-homogeneous curve on
algebraic restrictions is much more complicated.

{\bf Acknowledgements.} The author wishes to express his thanks to
M. Zhitomirskii for suggesting the subject and for many helpful
conversations and remarks during the writing of this paper. The
author thanks  Z. Jelonek for very useful remarks on the proof of
Theorem \ref{main-alg} and the referee of this paper for many
valuable suggestions.

\section{Quasi-homogeneity}
In this section we present the basic definitions and properties of
quasi-homogeneous germs.
\begin{defn}
 A curve-germ $f:(\mathbb
R,0)\rightarrow (\mathbb R^m,0)$ is {\bf quasi-homogeneous} if
there exist coordinate systems $t$ on $(\mathbb R,0)$ and
$(x_1,\cdots,x_m)$ on $(\mathbb R^m,0)$ and positive integers
$(\lambda_1,\cdots, \lambda_m)$ such that
$$df\left(t\frac{d}{dt}\right)=E\circ f,$$
where $E=\sum_{i=1}^m\lambda_i x_i\frac{\partial}{\partial x_i}$
is the germ of the Euler vector field on $(\mathbb R^m,0)$. The
coordinate system $(x_1,\cdots, x_m)$ is called quasi-homogeneous,
and numbers $(\lambda_1,\cdots, \lambda_m)$ are called weights.
\end{defn}

\begin{defn}
Positive integers $\lambda_1,\cdots,\lambda_m$ are {\bf linearly
dependent over non-negative integers} if there exists $j$ and
non-negative integers $k_i$ for $i\ne j$ such that
$\lambda_j=\sum_{i\ne j}k_i \lambda_i$. Otherwise we say that
$\lambda_1,\cdots,\lambda_m$ are {\bf linearly independent over
non-negative integers}.
\end{defn}

It is easy to see that quasi-homogeneous curves have the following
form in the quasi-homogeneous coordinates.
\begin{prop}
A curve-germ $f$ is quasi-homogeneous if and only if $f$ is
$\mathcal A$-equivalent to
$$t\mapsto (t^{\lambda_1},\cdots,t^{\lambda_k},0,\cdots,0),$$
where  $\lambda_1 < \cdots < \lambda_k$ are positive integers
linearly independent over non-negative integers.

$\lambda_1,\cdots \lambda_k$ generate the semigroup of the curve
$f$, which we denote by $(\lambda_1,\cdots, \lambda_k)$.
\end{prop}

The weights $\lambda_1,\cdots, \lambda_k$ are determined by $f$,
but weights $\lambda_{k+1},\cdots,\lambda_m$ can be arbitrary
positive integers. Actually in the next sections we study the
projection of $f$ to non zero components: $\mathbb R \ni t\mapsto
(t^{\lambda_1},\cdots, t^{\lambda_k})\in \mathbb R^{k}$.

\begin{defn}
The germ of a function, a differential $k$-form, or a vector field
$\alpha$ on $(\mathbb R^m,0)$ is {\bf quasi-homogeneous} in a
coordinate system $(x_1,\cdots,x_m)$ on $(\mathbb R^m,0)$ with
positive weights $(\lambda_1,\cdots, \lambda_m)$ if $\mathcal L_E
\alpha=\delta \alpha$, where $E=\sum_{i=1}^m\lambda_i
x_i\frac{\partial}{\partial x_i}$ is the germ of the Euler vector
field on $(\mathbb R^m,0)$ and $\delta$ is a real number called
the quasi-degree.
\end{defn}

 It is easy to show that $\alpha$ is quasi-homogeneous in a coordinate system
 $(x_1,\cdots,x_m)$ with weights $(\lambda_1,\cdots, \lambda_m)$ if and only if
 $F_t^{\ast}\alpha=t^{\delta}\alpha$, where $F_t(x_1,\cdots,x_m)=(t^{\lambda_1} x_1,\cdots, t^{\lambda_m}x_m)$.
 Then germs of quasi-homogeneous functions of
quasi-degree $\delta$ are germs of weighted homogeneous
polynomials of degree $\delta$. The coefficient
$f_{i_1,\cdots,i_k}$ of the quasi-homogeneous differential
$k$-form $\sum f_{i_1,\cdots,i_k} dx_{i_1}\wedge\cdots \wedge
dx_{i_k}$ of quasi-degree $\delta$ is a weighted homogeneous
polynomial of degree $\delta-\sum_{j=1}^k \lambda_{i_j}$. The
coefficient $f_i$  of the quasi-homogeneous vector field
$\sum_{i=1}^m f_i\frac{\partial}{\partial x_i}$ of quasi-degree
$\delta$ is a weighted homogeneous polynomial of degree $\delta+
\lambda_{i}$.

\begin{prop} \label{degree-Lie}If $X$ is the
germ of a quasi-homogeneous vector field of quasi-degree $i$ and
$\omega$ is the germ of a quasi-homogeneous differential form of
quasi-degree $j$ then $\mathcal L_X\omega$ is the germ of a
quasi-homogeneous differential form of quasi-degree $i+j$.
\end{prop}

\begin{proof}Since
$\mathcal L_EX=[E,X]=i X$ and $\mathcal L_E\omega=j \omega$, we
have
$$
\mathcal L_E(\mathcal L_X \omega)=\mathcal L_X(\mathcal L_E
\omega)+\mathcal L_{[E,X]}\omega=\mathcal L_X(j \omega)+\mathcal
L_{iX}\omega=j\mathcal L_X\omega+i\mathcal L_X\omega=(i+j)\mathcal
L_X\omega.
$$
It implies that $\mathcal L_X a$ is quasi-homogeneous of
quasi-degree $i+j$.
\end{proof}

\section{The method of algebraic restrictions}\label{def-alg-rest}
In this section we present basic facts on the method of algebraic
restrictions. The proofs of all results of this section can be
found in \cite{DJZ2}.

Given the germ of a smooth manifold $(M,p)$ denote by $\Lambda
^k(M)$ the space of all germs at $p$ of differential $k$-forms on
$M$. Given a curve-germ $f:(\mathbb R,0)\rightarrow (M,p)$
introduce the following subspaces of $\Lambda ^p(M)$:
$$\Lambda ^p_{Im f}(M) = \{\omega \in \Lambda ^p(M): \ \ \omega|_{f(t)}=0 \ \text {for any} \ t\in \mathbb R \};$$
$$\mathcal A^p_0(Im f, M) = \{\alpha  + d\beta : \ \ \alpha
\in \Lambda _{Im f}^p( M), \ \beta \in \Lambda _{Im
f}^{p-1}(M).\}$$ The relation $\omega|_{f(t)}=0$ means that the
$p$-form $\omega $ annihilates any $p$-tuple of vectors in
$T_{f(t)}M$, i.e. all coefficients of $\omega $ in some (and then
any) local coordinate system vanish at the point $f(t)$.
\begin{defn}
\label{main-def} The {\bf algebraic restriction} of $\omega $ to
a curve-germ $f:\mathbb R\rightarrow M$ is the equivalence class
of $\omega $ in $\Lambda ^p(M)$, where the equivalence is as
follows: $\omega $ is equivalent to $\widetilde \omega $ if
$\omega - \widetilde \omega \in \mathcal A^p_0(Im f, M)$.
\end{defn}

\begin{nt}The algebraic restriction of the germ of a form $\omega $ on $(M,p)$ to
a curve-germ $f$ will be denoted by $[\omega ]_f$. Writing
$[\omega ]_f=0$ (or saying that $\omega $ has zero algebraic
restriction to $f$) we mean that $[\omega ]_f = [0]_f$, i.e.
$\omega \in A^p_0(Im f, M)$.
\end{nt}

\begin{rem}\label{image}
If $g=f\circ \phi$ for a local diffeomorphism $\phi$ of $\mathbb
R$ then the algebraic restrictions $[\omega]_f$ and $[\omega]_g$
can be identified, because $Im f=Im g$.
\end{rem}

Let $(M,p)$ and $(\widetilde M,\tilde p)$ be germs of smooth
equal-dimensional manifolds. Let $f:(\mathbb R,0)\rightarrow
(M,p)$ be  a curve-germ in $(M,p)$.  Let $\tilde f:(\mathbb
R,0)\rightarrow (\widetilde M,\tilde p)$ be  a curve-germ in
$(\widetilde M,\tilde p)$. Let $\omega$ be the germ of a $k$-form
on $(M,p)$ and $\widetilde \omega$ be the germ of a $k$-form on
$(\widetilde M,\tilde p)$.

\begin{defn}Algebraic restrictions $[\omega ]_f$
and $[\widetilde \omega ]_{\widetilde f}$ are called {\bf
diffeomorphic} if there exists the germ of a diffeomorphism $\Phi:
(\widetilde M,\tilde p)\to (M,p)$ and the germ of a diffeomorphism
$\phi:(\mathbb R,0)\to (\mathbb R,0)$ such that $\Phi\circ \tilde
f\circ \phi=f$ and
$\Phi^{\ast}([\omega]_f):=[\Phi^{\ast}\omega]_{\Phi^{-1}\circ
f}=[\widetilde \omega]_{\tilde f}$.
\end{defn}

\smallskip

\begin{rem}
The above definition does not depend on the choice of $\omega$ and
$\widetilde \omega$ since a local diffeomorphism maps forms with
zero algebraic restriction to $f$ to forms with zero algebraic
restrictions to $\tilde f$.  If $(M,p)=(\widetilde M,\tilde p)$
and $f = \tilde f$ then the definition of diffeomorphic algebraic
restrictions reduces to the following one: two algebraic
restrictions $[\omega ]_f$ and $[\widetilde \omega ]_f$ are
diffeomorphic if there exist germs of diffeomorphisms $\Phi$ of
$(M,p)$ and $\phi$ of $(\mathbb R,0)$ such that $\Phi\circ f\circ
\phi=f$ and $[\Phi ^*\omega ]_f = [\widetilde \omega ]_f$.
\end{rem}

The method of algebraic restrictions applied to singular
quasi-homogeneous curves is based on the following theorem.
\begin{thm}[Theorem A in \cite{DJZ2}]
\label{DG}Let $f:(\mathbb R,0)\rightarrow (\mathbb R^{2n},0)$ be
the germ of a quasi-homogeneous curve. If $\omega _0, \omega _1$
are germs of symplectic forms on $(\mathbb R^{2n},0)$ with the
same algebraic restriction to $f$ then there exists the germ of a
diffeomorphism $\Phi:(\mathbb R^{2n},0)\rightarrow (\mathbb
R^{2n},0)$ such that $\Phi \circ f = f$ and $\Phi ^*\omega _1 =
\omega _0$.

 Two germs of
quasi-homogeneous curves $f, g$ of a fixed symplectic space
$(\mathbb R^{2n}, \omega)$ are symplectically equivalent if and
only if the algebraic restrictions of the symplectic form $\omega
$ to $f$ and $g$ are diffeomorphic.
\end{thm}

Theorem \ref{DG} reduces the problem of symplectic classification
of singular quasi-homogeneous curves to the problem of
diffeomorphic classification of algebraic restrictions of
symplectic forms to a singular quasi-homogeneous curve.

In the next section we prove that the set of algebraic
restrictions of 2-forms to a singular quasi-homogeneous curve is a
finite dimensional vector space. We now recall basic properties of
algebraic restrictions which are useful for a description of this
subset (\cite{DJZ2}).

Let $f$ be a quasi-homogeneous curve on $(\mathbb R^{2n},0)$.

First we can reduce the dimension of the manifold we consider due
to the following propositions.

\begin{prop}
\label{reduction} Let $(M,0)$ be the germ of a smooth submanifold
of $(\mathbb R^m,0)$ containing $Im f$. Let $\omega _1, \omega _2$
be germs of $k$-forms on $(\mathbb R^m,0)$. Then $[\omega _1]_f =
[\omega _2]_f$ if and only if $\bigl[\omega _1\vert _{TM}\bigr]_f
= \bigl[\omega _2 \vert _{TM}\bigr]_f$.
\end{prop}

\begin{prop}
\label{main-reduction} Let $f_1,f_2$ be curve-germs in $(\mathbb
R^m,0)$ whose images are contained in germs of equal-dimensional
smooth submanifolds $(M_1,0)$, $(M_2,0)$ respectively. Let $\omega
_1, \omega _2$ be germs of $k$-forms on $(\mathbb R^m,0)$. The
algebraic restrictions $[\omega _1]_{f_1}$ and $[\omega _2]_{f_2}$
are diffeomorphic if and only if the algebraic restrictions
$\bigl[\omega _1\vert _{TM_1}\bigr]_{f_1}$ and $\bigl[\omega
_2\vert _{TM_2}\bigr]_{f_2}$ are diffeomorphic.
\end{prop}

To calculate the space of algebraic restrictions of $2$-forms we
will use the following obvious properties.
\begin{prop}\label{d-wedge}
If $\omega\in \mathcal A_0^k(Im f,\mathbb R^{2n})$ then $d\omega
\in \mathcal A_0^{k+1}(Im f,\mathbb R^{2n})$ and $\omega\wedge
\alpha \in \mathcal A_0^{k+p}(Im f,\mathbb R^{2n})$ for any
$p$-form $\alpha$ on $\mathbb R^{2n}$.
\end{prop}

The next step of our calculation is the description of the
subspace of algebraic restriction of closed $2$-forms. The
following proposition is very useful for this step.
\begin{prop}
\label{th-all-closed} Let $a_1,\dots , a_k$ be a basis of the
space of algebraic restrictions of $2$-forms  to $f$ satisfying
the following conditions
\begin{enumerate}
\item $da_1 = \cdots = da_j = 0$, \item the algebraic restrictions
$da_{j+1}, \dots , da_k$ are linearly independent.
\end{enumerate}
 Then $a_1, \dots , a_j$ is a basis of the space
of algebraic restriction of closed $2$-forms to $f$.
\end{prop}

Then we need to determine which algebraic restrictions of closed
$2$-forms are realizable by symplectic forms. This is possible due
to the following fact.
\begin{prop}\label{rank}
Let $r$ be the minimal dimension of germs of smooth submanifolds
of $(\mathbb R^{2n},0)$ containing $Im f$. Let $(S,0)$ be one of
such germs of $r$-dimensional smooth submanifolds. Let $\theta$ be
the germ of a closed $2$-form on $(\mathbb R^{2n},0)$. There
exists the germ of a symplectic form $\omega$ on $(\mathbb
R^{2n},0)$ such that $[\theta]_f=[\omega]_f$ if and only if $rank
\theta|_{T_0 S}\ge 2r-2n$.
\end{prop}

\section{Discrete symplectic invariants.}\label{discrete}

Some new discrete symplectic invariants can be effectively
calculated using algebraic restrictions. The first one is a
symplectic multiplicity (\cite{DJZ2}) introduced  in \cite{IJ1} as
a symplectic defect of a curve $f$.

\begin{defn}
\label{def-mu}
 The {\bf symplectic multiplicity} $\mu_{sympl}(f)$ of a curve $f$ is the codimension of
 a symplectic orbit of $f$ in an $\mathcal A$-orbit of $f$.
\end{defn}

The second one is the index of isotropness \cite{DJZ2}.

\begin{defn}
The {\bf index of isotropness} $\iota(f)$ of $f$ is the maximal
order of vanishing of the $2$-forms $\omega \vert _{TM}$ over all
smooth submanifolds $M$ containing $Im f$.
\end{defn}

They can be described in terms of algebraic restrictions
\cite{DJZ2}.

\begin{prop}\label{sm}
The symplectic multiplicity  of a quasi-homogeneous curve $f$ in a
symplectic space is equal to the codimension of the orbit of the
algebraic restriction $[\omega ]_f$ with respect to the group of
local diffeomorphisms preserving $f$  in the space of algebraic
restrictions of closed  $2$-forms to $f$.
\end{prop}

\begin{prop}\label{ii}
The index of isotropness  of a quasi-homogeneous curve $f$ in a
symplectic space $(\mathbb R^{2n}, \omega )$ is equal to the
maximal order of vanishing of closed $2$-forms representing the
algebraic restriction $[\omega ]_f$.
\end{prop}

The above invariants are defined for the image of $f$. They have
the natural generalization to any subset of the symplectic space
\cite{DJZ2}.

There is one more discrete symplectic invariant introduced in
\cite{Ar1} which is defined specifically for a parameterized
curve. This is the maximal tangency order of a curve $f$ to a
smooth Lagrangian submanifold. If $H_1=...=H_n=0$ define a smooth
submanifold $L$ in the symplectic space then the tangency order of
a curve $f:\mathbb R\rightarrow M$ to $L$ is the minimum of the
orders of vanishing at $0$ of functions $H_1\circ f,\cdots,
H_n\circ f$. We denote the tangency order of $f$ to $L$ by
$t(f,L)$.

\begin{defn}
The {\bf Lagrangian tangency order} $Lt(f)$ of a curve $f$ is the
maximum of $t(f,L)$ over all smooth Lagrangian submanifolds $L$ of
the symplectic space.
\end{defn}

For a quasi-homogeneous curve $f$ with the semigroup
$(\lambda_1,\cdots,\lambda_k)$ the Lagrangian tangency order is
greater than $\lambda_1$.

$Lt(f)$ is related to the index of isotropness. If the index of
isotropness of $\omega$ to $f$ is $0$ then there does not exist a
closed $2$-form vanishing at $0$ representing algebraic
restriction of $\omega$. Then it is easy to see that the order of
tangency of $f$ to $L$ is not greater then $\lambda_k$.

The Lagrangian tangency order of a quasi-homogeneous curve in a
symplectic space can also be  expressed in terms of algebraic
restrictions.

The order of vanishing of the germ of a $1$-form $\alpha$ on a
curve-germ $f$ at $0$ is the minimum of the orders of vanishing of
functions $\alpha(X)\circ f$ at $0$ over all germs of smooth
vector fields $X$. If $\alpha=\sum_{i=1}^m g_i dx_i$ in local
coordinates $(x_1,\cdots,x_m)$ then the order of vanishing of
$\alpha$ on $f$ is the minimum of the orders of vanishing of
functions $g_i\circ f$ for $i=1,\cdots,m$.

\begin{prop}\label{lto}
Let $f$ be the germ of a quasi-homogeneous curve such that the
algebraic restriction of a symplectic form to it can be
represented by a closed $2$-form vanishing at $0$. Then the
Lagrangian tangency order of the germ of a quasi-homogeneous curve
$f$ is the maximum of the order of vanishing on $f$ over all
$1$-forms $\alpha$ such that $[\omega]_f=[d\alpha]_f$
\end{prop}

\begin{proof} Let $L$ be the germ of a smooth Lagrangian submanifold
in a standard symplectic space $(\mathbb
R^{2n},\omega_0=\sum_{i=1}^n dp_i\wedge dq_i)$. Then there exist
disjoint subsets $J,K\subset \{1,\cdots,n\}$, $J\cup
K=\{1,\cdots,n\}$ and a smooth function $S(p_J,q_K)$ (\cite{ArGi})
such that
\begin{equation}\label{L}
L=\{q_j=-\frac{\partial S}{\partial p_j}(p_J,q_K),
p_k=\frac{\partial S}{\partial q_k}(p_J,q_K), \ j\in J, \ k\in
K\}.
\end{equation}
It is obvious that the order of tangency of $f$ to $L$ is equal to
the order of vanishing of the following $1$-form:
$\alpha=\sum_{k\in K} p_k dq_k-\sum_{j\in J} q_j dp_j   -
dS(p_J,q_k)$ and $d\alpha=\omega_0$.

If closed $2$-forms have the same algebraic restrictions to $f$
then their difference can be written as a differential of a
$1$-form vanishing on $f$ by relative Poincare lemma for
quasi-homogeneous varieties \cite{DJZ1}. That implies that the
maximum of orders of vanishing of $1$-forms $\alpha$ on $f$
depends only on the algebraic restriction of $\omega=d\alpha$. Let
$f(t)=(t^{\lambda_1},\cdots,t^{\lambda_k},0,\cdots, 0)$. We may
assume that $[\omega]_f$ may be identified with $[d\alpha]_f$,
where $\alpha$ is a $1$-form on $\{x_{k+1}=\cdots =x_{2n}=0\}$ and
$d\alpha|_0=0$. In local coordinates $\alpha=\sum_{i=1}^k g_idx_i$
where $g_i$ are smooth function-germs. Let $\sigma$ be the
following germ of a symplectic form
$$
\sigma=d\alpha+\sum_{i=1}^k dx_i\wedge
dx_{k+i}+\sum_{i=1}^{n-k}dx_{2k+i}\wedge dx_{n+k+i}.
$$
Let $L$ be the following germ of a smooth Lagrangian submanifold
(with respect to $\sigma$)
$$
\{x_{k+i}=g_i,  \ i=1,\cdots,k, \ x_{2k+j}=0, \ j=1,\cdots, n-k\}.
$$
The tangency order of $f$ to $L$ is the same as the order of
vanishing of $\alpha$ on $f$.

It is obvious that the pullback of $\sigma$ to $\{x_{k+1}=\cdots
=x_{2n}=0\}$ is $d\alpha$. Then by Darboux-Givental theorem
(\cite{ArGi}) there  exists a local diffeomorphism which is the
identity on $\{x_{k+1}=\cdots =x_{2n}=0\}$ and maps $\sigma$ to
$\omega$. $L$ is mapped to a smooth Lagrangian submanifold (with
respect to the symplectic form $\omega$) with the same tangency
order to $f$.\end{proof}

\section{The proof of Theorem \ref{main-alg}}\label{proof}

In this section we prove Theorem \ref{main-alg}. The proof is
based on the following lemmas.

\begin{lem}\label{power}
Let $N$ be the germ of a subset of $\mathbb K^m$ at $0$. Let
$(x_1,\cdots,x_m)$ be a local coordinate system on $\mathbb K^m$.

The space of algebraic restrictions of $2$-forms to $N$ is finite
dimensional if and only if there exists a non-negative integer $L$
such that $x_i^L dx_j\wedge dx_k$ has zero algebraic restriction
to $N$ for any $i,j,k=1,\cdots,m$.
\end{lem}

\begin{proof}[Proof of Lemma \ref{power}] To prove the "only if"
part notice that  there exists an non-negative integer $K$ such
that algebraic restrictions
$$
[dx_j\wedge dx_k]_N, \ [x_i dx_j\wedge dx_k]_N, [x_i^2 dx_j\wedge
dx_k]_N,\cdots, \ [x_i^Kdx_j\wedge dx_k]_N
$$
are linearly dependent, since the space of algebraic restrictions
of $2$-forms to $N$ is finite dimensional. Therefore there exist a
non-negative integer $M$ and $c_1,\cdots,c_s \in \mathbb K$ such
that $[x_i^M(1+\sum_{l=1}^s c_sx_i^{l})dx_j\wedge dx_k]_N=0$. It
implies that $[x_i^Mdx_j\wedge dx_k]_N=0$. Now it is easy to see
that $L$ is the maximum of $M$ for all choices of $i,j,k$.

To prove the ''if'' part first notice that  any germ of a $2$-form
can be written in the local coordinates as $\sum_{1\le j<k\le
m}F_{j,k}(x)dx_j\wedge dx_k$, where $F_{j,k}(x)$ are
function-germs on $\mathbb K^m$. Using Taylor expansions of
$F_{j,k}(x)$ with the reminder of degree greater than $mL$ we
obtain the result, since $x_1^{i_1}\cdots x_m^{i_m} dx_j\wedge
dx_k$ has zero algebraic restriction to $N$ for $i_1+\cdots+i_m\ge
mL$.

\end{proof}

\begin{lem}\label{Wall}
Let $f:(\mathbb K,0)\rightarrow (\mathbb K^2,0)$ be the germ of a
$\mathbb K$-analytic parameterized curve in $\mathbb K^2$. Let
$(y,z)$ be a local coordinate system on $\mathbb K^2$, such that
the line $\{y=0\} $ does not contain $f(\mathbb K)$. Then there
exists a $\mathbb K$-analytic function-germ $H$ vanishing on $f$
of the following form $H(y,z)=z^p-G(y,z)y^l$, where $G$ is a
$\mathbb K$-analytic function-germ on $\mathbb K^2$, and $p, l$
are positive integers.
\end{lem}
\begin{proof}[Proof of Lemma \ref{Wall}] We use the method of a construction of $H$ described in
\cite{W} (the proof of Lemma 2.3.1 on page 28). $f$ is $\mathbb
K$-analytic then there exists a coordinate system $t$ on $\mathbb
K$  such that $f(t)=(t^m,\sum_{i=k}^{\infty} a_i t^i)$. We write
it in the following way $y=t^m, \ z=\sum_{i=k}^{\infty} a_i t^i$.
Any non-negative integer $i$ can be written in the following way
$i=qm+r$, where $r, q$ are integers such that $0\le r \le m-1$ and
$q \ge 0$. Thus $z=\sum_{r=0}^{m-1}t^r\phi_r(y)$, where
$\phi_r(y)=\sum_{q=0}^{\infty}a_{q m +r}y^q$ is $\mathbb
K$-analytic for $r=0,1,\cdots,m-1$. Then regard the equations
\begin{equation}\label{system}
t^az=\sum_{r=0}^{m-a-1}t^{a+r}\phi_r(y)+\sum_{r=m-a}^{m-1}t^{a+r-m}z\phi_r(y)
\ {\rm for} \ a=0,1,\cdots,m-1.
\end{equation}
 as a system of linear equations
for the unknowns $t^r$ $r=0,\cdots,m-1$ with coefficients in
$\mathbb K\{y,z\}$. The determinant $D(x,y)$ of this system has
the following form
$$\det \left[
\begin{array}{ccccccccccccccc}
z-\phi_0(y)&,& -\phi_1(y)&,&-\phi_2(y)&,& \cdots &,&-\phi_{m-1}(y)\\
-z\phi_{m-1}(y)&,&z-\phi_0(y)&,& -\phi_1(y)&,& \cdots &,&-\phi_{m-2}(y)\\
-z\phi_{m-2}(y)&,&-z\phi_{m-1}(y)&,&z-\phi_0(y)&,& \cdots &,&-\phi_{m-3}(y)\\
\vdots& &\vdots& &\vdots& &\ddots& &\vdots\\
-z\phi_1(y)&,&-z\phi_2(y)&,&-z\phi_3(y)&,& \cdots &,&z-\phi_0(y)\\
\end{array}\right]=
$$

$$
=z^m+\psi_1(y)z^{m-1}+\cdots+\psi_{m-1}(y)z+\psi_m(y),
$$
where $\psi_1,\cdots,\psi_m$ are $\mathbb K$-analytic
function-germs. Since the values $t^r$ for $r=0,\cdots,m-1$
provide non-zero solutions of (\ref{system}), the determinant
$D(y,z)$ vanishes on the image of the curve $f$. Since $f(0)=0$ we
have that $\psi_m(0)=0$.

Thus we can decompose $D(y,z)$ in the following form
$$
D(y,z)=z^m+\psi_1(0)z^{m-1}+\cdots+\psi_{m-1}(0)z+y^lF(y,z)=z^ph(z)+y^lF(y,z),
$$
where $h$ is a polynomial of degree $m-k$ that does not vanish at
$0$, $p, l$ are positive integers and $F$ is a $\mathbb
K$-analytic function-germ. Now we take
$H(y,z)=\frac{D(y,z)}{h(z)}$.
\end{proof}

\begin{lem} \label{Jelonek} Let $C$ be the  germ of a $\mathbb K$-analytic curve on
$\mathbb K^2$ at $0$. Let $(y,z)$ be a local coordinate system on
$\mathbb K^2$, such that the line $\{y=0\} $ does not contain any
branch of  $C$. Then there exists a $\mathbb K$-analytic
function-germ $H$ vanishing on $f$ of the following form
$H(y,z)=z^p-G(y,z)y^l$, where $G$ is a $\mathbb K$-analytic
function-germ on $\mathbb K^2$ and $p,l$ are positive integers.
\end{lem}

\begin{proof}We decompose $C$ into branches $C_1,\cdots, C_s$. Then we apply Lemma
\ref{Wall} to each branch $C_i$. We obtain  a $\mathbb K$-analytic
function-germ vanishing on $C_i$ of the form
$H_i(y,z)=z^{p_i}-G_i(y,z)y^{l_i}$, where $p_i, l_i$ are positive
integers and $G_i$ is  a $\mathbb K$-analytic function-germ for
$i=1,\cdots,s$. Now we may take $H=H_1\cdots H_s$, which vanishes
on $C$ and has the desired form.
\end{proof}
\begin{lem}\label{Tjurina} Let $N$ be the germ of a subset of $\mathbb K^2$ at $0$. Let $H$ be a
$\mathbb K$-analytic  function-germ on $\mathbb K^2$ vanishing on
$N$.

If $H$ has a regular point at $0$ or an isolated critical point at
$0$ then the space of algebraic restrictions of $2$-forms to $N$
is finite dimensional.
\end{lem}

\begin{proof} The space of algebraic restrictions of $2$-forms to $\{H=0\}$ is isomorphic to $\frac{\mathcal C_2}{<H,\nabla
H>}$ \cite{DJZ2}, where $\mathcal C_2$ is the space of $\mathbb
K$-analytic function-germs on $\mathbb K^2$. Thus its dimension is
finite and equal to the Tjurina number of $\{H=0\}$. $N$ is a
subset of $\{H=0\}$. Hence the dimension of the space of algebraic
restriction of $2$-forms to $N$ is smaller than the Tjurina number
of $\{H=0\}$, and consequently it is finite.
\end{proof}

\begin{proof}[Proof of Theorem \ref{main-alg}] Let $C$ be the germ
of a $\mathbb K$-analytic curve in $\mathbb K^m$ at $0$. In fact
we prove that the vector space of algebraic restrictions of all
$2$-forms to $C$ is finite dimensional. It is obvious that the set
of algebraic restrictions of closed $2$-forms is a vector subspace
of the vector space of algebraic restrictions of all $2$-forms.

Let $(x_1,\cdots,x_m)$ be a coordinate system on $\mathbb K^m$ and
let
 $$\pi_{j,k}:\mathbb K^m\ni (x_1,\cdots,x_m) \rightarrow
(x_j,x_k)\in \mathbb K^2$$ be the standard projection. We choose a
coordinate system in such way that for any $j\ne k$ $\pi_{j,k}(C)$
is the germ of a $\mathbb K$-analytic curve on $\mathbb K^2$ at
$0$ such that the lines $\{x_j=0\}$ and $\{x_k=0\}$ do not contain
any branch of $\pi_{j,k}(C)$.

Then the space of algebraic restrictions of $2$-forms to
$\pi_{j,k}(C)$ is finite dimensional by Lemma \ref{Tjurina}, since
$\pi_{j,k}(C)$ may have a non-singular point at $0$ or an isolated
singular point at $0$. By Lemma \ref{power} there exists a
positive integer $K$ such that $x_j^Kdx_j\wedge dx_k$ has zero
algebraic restriction to $\pi_{j,k}(C)$ and consequently   it has
zero algebraic restriction to $C$.

 By Lemma
\ref{Jelonek} there exist positive integers $p,l$ and a $\mathbb
K$-analytic function-germ $G$ on $\mathbb K^2$  such that the
function-germ  $H(x_j,x_i)=x_i^p+G(x_j,x_i)x_j^l$ vanishes on
$\pi_{j,i}(C)$ and  consequently it vanishes on $C$.  It implies
that
$$
x_i^{pK}dx_j\wedge dx_k=(-G(x_j,x_i))^Kx_j^{lK}dx_j\wedge dx_k
$$
 has
zero algebraic restriction to $C$ too.

Hence by Lemma \ref{power} we obtain that the space of algebraic
restrictions of $2$-forms to $C$ is finite dimensional.
\end{proof}

\section{Quasi-homogeneous algebraic restrictions}\label{qh-alg}
In this section we prove  that the action by diffeomorphisms
preserving the curve is totaly determined by infinitesimal action
by liftable vector fields and the space of such vector fields
which act nontrivially on algebraic restrictions is a finite
dimensional vector space spanned by quasi-homogeneous liftable
vector fields of bounded quasi-degrees.

The proof of Theorem \ref{main-alg} is very easy in the case of
quasi-homogeneous parameterized curves. Let $f$ be the germ of a
quasi-homogeneous curve. Then $f$ is $\mathcal A$-equivalent to
$f(t)=(t^{\lambda_1},\cdots,t^{\lambda_k},0,\cdots,0)$. By
Proposition \ref{reduction} we consider  forms in $x_1,\cdots,x_k$
coordinates only . We may also assume that the greatest common
divisor $g(\lambda_1,\cdots, \lambda_k)$ is $1$. If it is not $1$
we introduce weights $\lambda_i/gcd(\lambda_1,\cdots,\lambda_k)$
for $x_i$, $i=1,\cdots,k$. The proof of Theorem \ref{main-alg} in
this special case is based on the following easy observation.
\begin{lem} \label{li-lj}The function-germ
$h(x)=x_i^{\lambda_j}-x_j^{\lambda_i}$ vanishes on $f$.
\end{lem}
The above lemma implies the following facts.
\begin{lem}\label{lilj}
The $2$-form $x_i^{\lambda_j-1}dx_i\wedge dx_j$ has zero
algebraic restriction to $f$
\end{lem}
\begin{proof}[Proof of Lemma \ref{lilj}]  By Lemma \ref{li-lj} $dh$
has zero algebraic restriction to $f$. It implies that
$\frac{1}{\lambda_j}dh\wedge dx_j=x_i^{\lambda_j-1}dx_i\wedge
dx_j$ has zero algebraic restriction to $f$.
\end{proof}

\begin{lem} \label{dx_idx_j}
If the monomials $s(x)=\prod_{l=1}^k x_l^{s_l}$ and
$p(x)=\prod_{l=1}^k x_l^{p_l}$ have the same quasi-degree then the
forms $s(x)dx_i\wedge dx_j$ and $p(x)dx_i\wedge dx_j$ have the
same algebraic restrictions to $f$.
\end{lem}
\begin{proof}[Proof of Lemma \ref{dx_idx_j}] The function-germ $s(x)-p(x)$ vanishes on $f$.
\end{proof}

The above lemmas imply that we can choose the quasi-homogeneous
bases of the space of algebraic restrictions of $2$-forms to $f$.
Thus as a corollary of Theorem \ref{main-alg} and the above lemmas
we obtain the following theorem.

\begin{thm}\label{finite}
The space of algebraic restrictions of closed $2$-forms to the
germ of a quasi-homogeneous curve $f$ is a finite dimensional
vector space spanned by algebraic restrictions of
quasi-homogeneous closed $2$-forms of bounded quasi-degrees.
\end{thm}

We will use quasi-homogeneous grading on the space of algebraic
restrictions. Therefore we define quasi-homogeneous algebraic
restrictions.

Let $f$ be the germ of  a quasi-homogeneous curve on $(\mathbb
R^m,0)$. Let $\omega$ be the germ of a $k$-form  on $(\mathbb
R^m,0)$. By $\omega^{(r)}$ we denote a quasi-homogeneous part of
quasi-degree $r$ in the Taylor series of $\omega$.  It is easy to
see that if $h$ is a function-germ  on $(\mathbb R^m,0)$ and
$h\circ f=0$ then $h^{ (r)}\circ f=0$ for any $r$.  This simple
observation implies the following proposition.
\begin{prop}\label{qh-al-rest}
If $[\omega]_f=0$ then $[\omega^{(r)}]_f=0$ for any $r$.
\end{prop}
Proposition \ref{qh-al-rest} allows to define quasi-homogeneous
algebraic restriction.

\begin{defn}
Let $a=[\omega]_f$ be an algebraic restriction to $f$. The
algebraic restriction $a^{(r)}=[\omega^{[r]}]_f$ is called the
{\bf quasi-homogeneous part of quasi-degree} $r$ of an algebraic
restriction $a$. $a$ is {\bf quasi-homogeneous} of quasi-degree
$r$ if $a=a^{(r)}$.
\end{defn}

We consider the action on the space of algebraic restrictions of
closed $2$-forms by the group of diffeomorphism-germs which
preserve the curve $f$ to obtain a complete symplectic
classification of curves (Theorem \ref{DG}). The tangent space at
the identity to this group is given by the space of  vector fields
liftable over $f$.

\begin{defn}[\cite{Za}, \cite{BPW}]
The germ $X$ of a vector field on $(\mathbb R^m,0)$ is called {\bf
liftable} over $f$ if there exists a function germ $g$ on
$(\mathbb R,0)$ such that
$$
g\left(\frac{df}{dt}\right)=X\circ f.
$$
\end{defn}

The tangent space to the orbit of an algebraic restriction $a$ is
given by $\mathcal L_Xa$ for all vector field $X$ liftable over
$f$.  The Lie derivative of an algebraic restriction with respect
to a liftable vector field is well defined due to the following
proposition.

\begin{prop}\label{zero-Lie-zero} Let $X$ be the germ of a  vector field on
$(\mathbb R^m,0)$ liftable over $f$ and $\omega$ be the germ of a
$k$-form on $(\mathbb R^m,0)$. If $[\omega]_f=0$ then $[\mathcal
L_X\omega]_f=0$.
\end{prop}
\begin{proof}
This is a consequence of the Cartan formula and the following
fact: $dh(X)\circ f=0$ for any function-germ $h$ on $(\mathbb
R^m,0)$ vanishing on $f$. To prove the above fact let us notice
that $dh(X)\circ f=(dh\circ f)(X\circ f)=(dh\circ
f)df\left(g\frac{d}{dt}\right)=d(h\circ
f)\left(g\frac{d}{dt}\right)$.
\end{proof}
By the Cartan formula we also obtain the following proposition.
\begin{prop}\label{0image}
If $X$ is a vector field vanishing on the image of $f$ then
$\mathcal L_Xa=0$ for any algebraic restriction $a$ to $f$.
\end{prop}
If $f$ is quasi-homogeneous then the Euler vector field $E$ is
liftable over $f$. The following proposition describes its
infinitesimal action on quasi-homogeneous algebraic restrictions.
\begin{prop}
If an algebraic restriction $a$ to $f$ is quasi-homogeneous of
quasi-degree $\delta$ then $\mathcal L_Ea=\delta a$.
\end{prop}

Let $X$ be a smooth vector field. By $X^{(r)}$ we denote the
quasi-homogeneous part of quasi-degree $r$ in the Taylor series of
$X$. We have the following result.

\begin{prop} If $X$ is liftable over $f$ then $X^{(r)}$ is
liftable over $f$.
\end{prop}
\begin{proof} We assume that
$f(t)=(t^{\lambda_1},\cdots,t^{\lambda_k},0,\cdots,0)$. Then
$X\circ f=g(t)df/dt$ for some function-germ $g$ on $\mathbb R$. It
implies that
$$X^{(r)}\circ
f=\frac{1}{(r+1)!}\frac{d^{r+1}g}{dt^{r+1}}(0)t^{r+1}\frac{df}{dt}.
$$
\end{proof}

Let $K(f)$ be the minimal natural number such that all
quasi-homogeneous algebraic restrictions to $f$ of closed
$2$-forms of quasi-degree greater than $K(f)$ vanish. By Theorem
\ref{finite} $K(f)$ is finite.

\begin{thm}\label{vf-finite}
Let $f(t)=(t^{\lambda_1},\cdots,t^{\lambda_k},0,\cdots,0)$. Let
$X_s$ be the  germ of a vector field such that $X_s\circ
f=t^{s+1}df/dt$. Then the tangent space to the orbit of the
quasi-homogeneous algebraic restriction $a_r$ of quasi-degree $r$
is spanned by $L_{X_s}a_r$ for $s$ that are  $\mathbb Z_{\ge
0}$-linear combinations of $\lambda_1,\cdots,\lambda_k$ and
smaller than $K(f)-r$.
\end{thm}
\begin{proof}
Let $Y$ be a quasi-homogeneous vector field liftable over $f$.
Then $Y\circ f=c t^{s+1} df/dt$ where $s$ is the quasi-degree of
$Y$ and $c\in \mathbb R$. By Proposition \ref{0image} we obtain
that $L_Ya_r=cL_{X_s}a_r$, since $(Y-cX_s)\circ f=0$. If $Z$ is a
liftable vector field we can decompose it to
$\sum_{s=0}^{K(f)-r-1}Z^{(s)}+V$, where $V$ is a liftable vector
field such that $V^{(s)}=0$ for $s<K(f)-r$. Then
$L_Za_r=\sum_{s=0}^{K(f)-r-1}c_s L_{X_s}a_r + L_Va_r$, where
$c_s\in \mathbb R$ for $s=0,\cdots, K(f)-r-1$. Proposition
\ref{degree-Lie} implies that $(L_Va_r)^{(s)}=0$ for $s<K(f)$. By
Taylor expansion we can decompose $L_Va_r$ to $\sum_{i=1}^m
f_ib_i$, where $f_i$ are function-germs and $b_i$ are
quasi-homogeneous algebraic restrictions of quasi-degree greater
than $K(f)-1$. Thus $L_Va_r=0$.
\end{proof}

Theorem \ref{vf-finite} implies Theorem \ref{main-lift}. Now we
prove the following theorem which is crucial for the description
of the action of diffeomorphisms preserving $f$ on the space of
algebraic restrictions to $f$.
\begin{thm}\label{action}
Let $a_1, \cdots, a_p$ be a quasi-homogeneous basis of
quasi-degrees
$\delta_1\le\cdots\le\delta_s<\delta_{s+1}\le\cdots\le\delta_p$ of
the space of algebraic restrictions of closed $2$-forms to $f$ .
Let $a=\sum_{j=s}^p c_ja_j$, where $c_j\in \mathbb R$ for
$j=s,\cdots, p$ and $c_s\ne 0$.

If there exists a liftable quasi-homogeneous vector field $X$ over
$f$ such that $\mathcal L_Xa_s=ra_k$ for $k>s$ and $r\ne 0$ then
$a$ is diffeomorphic to $\sum_{j=s}^{k-1} c_ja_j+\sum_{j=k+1}^p
b_j a_j$, for some $b_j\in \mathbb R, \ j=k+1,\cdots,p$.
\end{thm}
\begin{proof}We use the Moser homotopy method. Let
$a_t=\sum_{j=s}^{k-1}
c_ja_j+(1-t)c_ka_k+\sum_{j=k+1}^{s}b_j(t)a_j$ where $b_j(t)$ are
smooth functions $b_j:[0;1]\rightarrow \mathbb R$ such that
$b_j(0)=c_j$ for $j=k+1,\cdots,p$. Let $\Phi_t, \ t\in[0;1]$ be a
flow of the vector field $\frac{c_k}{rc_s}V$ . We show that there
exist such functions $b_j$ that
\begin{equation}\label{flow}
\Phi_t^{\ast}a_t=a
\end{equation} for $t\in [0;1]$. Differentiating (\ref{flow}) we obtain
\begin{equation}\label{ODE}
\mathcal L_{\frac{c_k}{rc_s}X}a_t=c_ka_k-\sum_{j=k+1}^p
\frac{db_j}{dt}a_j.
\end{equation}
Since $\mathcal L_{X}a_s=ra_k$, the quasi-degree of $X$ is
$\delta_k-\delta_s$. Hence the quasi-degree of $\mathcal
L_{\frac{c_k}{rc_s}X}a_j$ is greater than $\delta_k$ for $j>s$.
Then $b_j$ are solutions of the system of $p-k$ first order linear
ODEs defined by (\ref{ODE}) with the initial data $b_j(0)=c_j$ for
$j=k+1,\cdots,p$. It implies that $a_0=a$ and
$a_1=\sum_{j=s}^{k-1} c_ja_j+\sum_{j=k+1}^{s}b_j(1)a_j$ are
diffeomorphic.
\end{proof}
\subsection{Remarks on the algorithm for a quasi-homogeneous parameterized curve with an arbitrary
semi-group}

The results of section \ref{qh-alg} allows us to formulate  an
algorithm for the classification of the symplectic singularities
of an arbitrary quasi-homogeneous parameterized curve-germ $f$. It
is obvious that this algorithm depends only on the semigroup of
the curve singularity.

Let us assume that the semigroup have the following form:
$$(\lambda_1,\cdots,\lambda_k),$$
where  $\lambda_1 < \cdots < \lambda_k$ are positive integers
linearly independent over non-negative integers. We use the
quasi-homogeneous grading on the space of algebraic restrictions
of $2$-forms with weights $(\lambda_1,\cdots,\lambda_k)$ for
coordinates $(x_1,\cdots,x_k)$. We may also assume that
$\lambda_1,\cdots, \lambda_k$ are relatively prime. If they are
not we introduce weights $\lambda_i/{\rm
gcd}(\lambda_1,\cdots,\lambda_k)$ for $x_i$, $i=1,\cdots,k$, where
${\rm gcd}(\lambda_1,\cdots,\lambda_k)$ is the greatest common
divisor of $\lambda_1,\cdots, \lambda_k$.

We fixed the quasi-degree $\delta$ starting with
$\lambda_1+\lambda_2$ since there are no quasi-homogeneous
$2$-forms with a smaller quasi-degree.

$2$-forms of the quasi-degree $\delta$ (together with the zero
$2$-form) form a finite dimensional vector subspace of the space
of differential $2$-forms.

By Lemma \ref{dx_idx_j} algebraic restrictions of all forms of the
quasi-degree $\delta$ of the form
\begin{equation}\label{gen-ij}
s(x)dx_i\wedge dx_j
\end{equation} for fixed $i\ne j$ are linearly
dependent. Hence for all $i<j$ we need to check if the equation
\begin{equation}\label{eq-ij}
a_1 \lambda_1+\cdots+a_k \lambda_k=d-\lambda_i-\lambda_j
\end{equation}
has a solution $a_1,\cdots, a_k$ in non-negative integers.

For fixed $i<j$ we take one of the solutions of equation
(\ref{eq-ij}) (if it exists). All other algebraic restrictions of
the $2$-forms of the form (\ref{gen-ij}) are linear combinations
of $[x_1^{a_1}\cdots x_k^{a_k}dx_i\wedge dx_j]_f$

To find a basis of algebraic restrictions of quasi-homogeneous
$2$-forms with the quasi-degree $\delta$ we are looking for
quasi-homogeneous functions vanishing on $f$ of a quasi-degree
$\delta$.

To find them we need to find all solutions of the equation
\begin{equation}\label{eq-i}
a_1 \lambda_1+\cdots+a_k \lambda_k=d-\lambda_i
\end{equation}
If $(a_1,\cdots,a_k)$ and $(b_1,\cdots,b_k)$ are distinct solution
of (\ref{eq-i}) then a function-germ
\begin{equation}\label{vanish-i} H(x_1,\cdots,
x_k)=x_1^{a_1}\cdots x_k^{a_k}-x_1^{b_1}\cdots x_k^{b_k}
\end{equation}
vanishes on $f$ and the form $dH\wedge dx_i$ has zero algebraic
restriction to $f$. In this way we obtain all relations between
algebraic restrictions of quasi-homogeneous forms of quasi-degree
$\delta$ and consequently we find a basis of this vector space.

Then we proceed to algebraic restrictions with quasi-degree
$\delta+1$.

For some quasi-degrees  we obtain that all quasi-homogeneous
$2$-forms have zero algebraic restriction to $f$. Then using the
fact that quasi-homogeneous forms of a sufficiently high
quasi-degree can be obtained by multiplication by functions of
quasi-homogeneous forms of lower degrees we get that all $2$-forms
of a sufficiently high quasi-degree have zero algebraic
restriction. In this way we construct the quasi-homogeneous basis
of the space of algebraic restriction of all $2$-forms.

Then by Proposition \ref{th-all-closed} we get the
quasi-homogeneous basis of the space of algebraic restriction of
closed $2$-forms from the quasi-homogeneous basis of the space of
algebraic restriction of all $2$-forms.

Then we calculate the number $K(f)$ and we find germs of vector
field such that $X_s\circ f=t^{s+1}df/dt$ for $s$ that are
representable as a non-negative integers combinations of
$\lambda_1,\cdots,\lambda_k$ and smaller than
$K(f)-\lambda_1-\lambda_2$. By Theorem \ref{vf-finite} the tangent
space to the orbit of the quasi-homogeneous algebraic restriction
$a_r$ of quasi-degree $r$ is spanned by $L_{X_s}a_r$.

Finally we apply Theorem \ref{action} to get the classification of
algebraic restrictions. From this classification we easily obtain
the symplectic singularities normal forms.

In the next sections we apply the above algorithm for curves with
semigroups $(3,4,5)$, $(3,5,7)$ and $(3,7,8)$.

Although the algorithm works very well for concrete examples, the
problem of calculations of the dimension of the space of algebraic
restrictions of closed $2$-forms to a quasi-homogeneous
parameterized curve in terms of the semigroup of this curve is
complicated. It is obvious that it is related to the classical
Frobenius coin problem (the diophantine Frobenius problem)
\cite{RA}.

\begin{problem}
Given $k$ relatively prime positive integers $\lambda_1,\cdots,
\lambda_k$, find the largest natural number (denoted by  ${\rm
g}(\lambda_1,\cdots, \lambda_k)$ and called {\bf Frobenius
number}) that is not representable as a non-negative integer
combination of $\lambda_1,\cdots, \lambda_k$.
\end{problem}

By Schur's theorem Frobenius number is finite (\cite{RA}). The
formula for Frobenius number for $k=2$ was found by J. J.
Sylvester: $g(\lambda_1,\lambda_2)=\lambda_1 \lambda_2 - \lambda_1
- \lambda_2$ (\cite{S}). Sylvester also demonstrated that there
are $(\lambda_1-1)(\lambda_2-1)/2$ non-representable natural
numbers. More complicated formulas and fast algorithms to
calculate Frobenius number for $k=3$ are known, but the general
problem for arbitrary $k$ is known to be NP-hard (\cite{RA}).

%
%
\section{Symplectic singularities of curves with the semigroup $(3,4,5)$} \label{345}
In this section we apply the results of the previous section to
prove the following theorem.

\begin{thm}\label{th(3,4,5)}Let $(\mathbb
R^{2n},\omega_0=\sum_{i=1}^n dp_i\wedge dq_i$) be the symplectic
space  with the canonical coordinates $(p_1,q_1,\cdots, p_n,q_n)$.

Then the germ of a curve $f:(\mathbb R,0)\rightarrow (\mathbb
R^{2n},0)$ with the semigroup $(3,4,5)$ is symplectically
equivalent to one and  only one of the curves presented in the
second column of the Table \ref{ss345} (on page 12) for $n>2$ and
$f$ is symplectically equivalent to one and only one of the curves
presented in the second column and rows 1-3 for $n=2$.

The symplectic multiplicity, the index of isotropness and the
Lagrangian tangency order are presented in the third, fourth and
fifth columns of Table \ref{ss345}.

\begin{table}[h]
\begin{center}
\begin{tabular}{|c|l|c|c|c|}

 \hline

   & normal form of $f$ & $\mu_{sympl}(f)$& $\iota(f)$ & $Lt(f)$\\ \hline

   1 &$t\mapsto(t^3,t^4,t^5,0,\cdots,0)$ & $0$ & $0$& $4$\\ \hline

 2 &$t\mapsto(t^3,\pm t^5,t^4,0,\cdots,0)$ &$1$  &$0$  &$5$\\ \hline

 3 & $t\mapsto(t^3,0,t^4,t^5,0,\cdots,0)$ & $2$ &$0$  &$5$\\ \hline

4 &$t\mapsto(t^3,\pm t^7,t^4,0,t^5,0,\cdots,0)$ &$3$  &$1$  &$7$\\
\hline

5 &$t\mapsto(t^3,t^8,t^4,0,t^5,0,\cdots,0)$ &$4$  &$1$  &$8$\\
\hline

6 & $t\mapsto(t^3,0,t^4,0,t^5,0,\cdots,0)$ & $5$ &$\infty$  &$\infty$\\
\hline

\end{tabular}
\end{center}
\caption{Symplectic classification of curves with the semigroup
$(3,4,5)$.}\label{ss345}
\end{table}

\end{thm}

We use the method of algebraic restrictions.
 The germ of a curve
$f:\mathbb R\ni t\mapsto f(t)\in \mathbb R^{2n}$ with the
semigroup $(3,4,5)$ is $\mathcal A$-equivalent to $t\mapsto
(t^3,t^4,t^5,0,\cdots,0)$. First we calculate the space of
algebraic restrictions of $2$-forms to the image of $f$ in
$\mathbb R^{2n}$.

\begin{prop}\label{basis-all}
The space of algebraic restrictions of differential $2$-forms to
$f$ is the $6$-dimensional vector space spanned by the following
algebraic restrictions:
$$
a_7=[dx_1\wedge dx_2]_g, \ a_8=[dx_3\wedge dx_1]_g, \
a_9=[dx_2\wedge dx_3]_g,
$$
$$
a_{10}=[x_1dx_1\wedge dx_2]_g, \ a_{11}=[x_2 dx_1\wedge dx_2]_g, \
a_{12}=[x_1dx_2\wedge dx_3]_g,
$$
where $\delta$ is quasi-degree of $a_{\delta}$.
\end{prop}

\begin{proof}The image of $f$ is contained in the following $3$-dimensional smooth
submanifold $\{x_{\ge 4}=0\}$. By Proposition \ref{reduction} we
can restrict our consideration to $\mathbb R^3$ and the curve
$g:\mathbb R\ni t \mapsto (t^3,t^4,t^5)\in \mathbb R^{3}$. $g$ is
quasi-homogeneous with weights $3, 4, 5$ for variables $x_1, x_2,
x_3$. We use the quasi-homogeneous grading on the space of
algebraic restrictions of differential $2$-forms to
$g(t)=(t^3,t^4,t^5)$ with these weights. It is easy to see that
the quasi-homogeneous functions or $2$-forms of a fixed
quasi-degree form a finite dimensional vector space. The same is
true for quasi-homogeneous algebraic restrictions of $2$-forms of
a fixed quasi-degree.

There are no quasi-homogeneous function-germs on $\mathbb R^3$
vanishing on $g$ of quasi-degree less than $8$. The vector space
of quasi-homogeneous function-germs of degree $i=8,9,10$ vanishing
on $g$ is spanned by $f_i$ presented in Table \ref{v345} together
with their differentials. We do not need to consider
quasi-homogeneous function-germs of higher quasi-degree, since
using $f_8$, $f_9$ and $f_{10}$ we show that algebraic
restrictions of quasi-homogeneous $2$-forms of quasi-degree
greater than $12$ are zero (see Table \ref{b345}) and all possible
relations of algebraic restrictions of quasi-homogeneous $2$-forms
of quasi-degree less than $13$ are generated by quasi-homogeneous
functions vanishing on $g$ of quasi-degree less than $13-3=10$.
\begin{table}
\begin{center}
\begin{tabular}{|c|c|c|}

 \hline

  quasi-degree $\delta$ & $f_{\delta}$& differential $df_{\delta}$\\ \hline

   $8$ & $x_1x_3-x_2^2$ & $x_1dx_3+x_3dx_1-2x_2dx_2$\\ \hline

  $9$ & $x_2x_3-x_1^3$ & $x_2dx_3+x_3dx_2-3x_1^2dx_1$\\ \hline

$10$ & $x_1^2x_2-x_3^2$ & $x_1^2dx_2+2x_1x_2dx_1-2x_3dx_3$\\
\hline
\end{tabular}
\end{center}
\caption{Quasi-homogeneous function-germs of quasi-degree $8, \ 9,
\ 10$ vanishing on the curve $t\mapsto
(t^3,t^4,t^5)$.}\label{v345}
\end{table}
Now we can calculate the space of algebraic restrictions of
$2$-forms. The scheme of the proof is presented in Table
\ref{b345}. The first column of this table contains possible
degree $\delta$ of a $2$-form. In the second column there is a
basis of the algebraic restrictions of $2$-forms of degree
$\delta$. In the third column we present the basis of $2$-forms of
degree $\delta$. In the fourth column we show the relations
between algebraic restrictions of elements of the basis of
$2$-forms of degree $\delta$. The last column contains the
sketches of proofs of these relations.

\begin{table}
\begin{center}
\begin{tabular}{|c|c|c|c|c|}

 \hline

  $\delta$&  basis & forms & relations & proof\\ \hline

   $7$ & $a_7$&$\alpha_7=dx_1\wedge dx_2$&$a_7:=[\alpha_7]_g$ &\\ \hline

  $8$ & $a_8$ & $\alpha_8=dx_3\wedge dx_1$&$a_8:=[\alpha_8]_g$ &\\ \hline

$9$ &  $a_9$& $\alpha_9=dx_2\wedge dx_3$ &$a_9:=[\alpha_9]_g$  &\\
\hline

$10$ &  $a_{10}$ & $x_1\alpha_7$& $a_{10}:=x_1a_7$&\\
\hline

$11$ & $a_{11}$ & $x_2\alpha_7, $ & $a_{11}:=x_2a_7$ &
$[df_8\wedge dx_1]_g=0$ \\
& & $x_1\alpha_8$&$a_{11}=-2x_1a_8$ &\\ \hline

$12$ & $a_{12}$ &$x_3\alpha_7,$ &
$a_{12}:=x_3a_7$   & $[df_9\wedge dx_1]_g=0$\\
& &$x_2\alpha_8, $ &$a_{12}=x_2a_8$ & $[df_8\wedge
dx_2]_g=0$ \\
& &$x_1 \alpha_9$ &$a_{12}=x_1a_9$ &\\ \hline

$13$ & $0$ &$x_1^2\alpha_7, $  &$x_1^2a_7=0$
  & $[df_{10}\wedge dx_1]_g=0$\\
& &$x_3\alpha_8,$ & $x_3a_8=0$& $[df_9\wedge dx_2]_g=0$ \\
& &$x_2\alpha_9$ & $x_2a_9=0$& $[df_8\wedge dx_3]_g=0$ \\
\hline

$14$ & $0$ &$x_1x_2\alpha_7, $  &$x_1x_2a_7=0$
  & $[df_{10}\wedge dx_2]_g=0$\\
& &$x_1^2\alpha_8,$ & $x_1^2a_8=0$& $[df_9\wedge dx_3]_g=0$ \\
& &$x_2\alpha_{10}$ & $x_2a_{10}=0$& $x_1[df_8\wedge dx_1]_g=0$ \\
\hline

$15$ & $0$ &$x_1x_3\alpha_7, $  &$x_1x_3a_7=0$ & $[df_{10}\wedge dx_3]_g=0$\\
& &$x_1x_2\alpha_8,$ & $x_1x_2a_8=0$& $x_1[df_9\wedge dx_1]_g=0$ \\
& &$x_1^2\alpha_9$ & $x_1^2a_9=0$& $x_1[df_8\wedge dx_2]_g=0$ \\
\hline

$\ge 16$ & $0$ &$x_1\beta_{\ge 13},$ &$b_{\ge 13}:=[\beta_{\ge
13}]_g$
  & $x_2b_{\ge12}=x_1b'_{\ge 13}$ \\
& & &$x_1b_{\ge 13}=0$&$x_3b_{\ge11}=x_1b''_{\ge 13} $ \\
& &$x_2\beta_{\ge 12},$ &$b_{\ge 12}:=[\beta_{\ge 12}]_g$ & $\delta(b_{\ge 13})\ge 13$ \\
& & &$x_2b_{\ge12}=0$ & $b_{\ge 13}=0$ \\
& &$x_3\beta_{\ge11}$ & $b_{\ge11}:=[\beta_{\ge11}]_g$ &  \\
& & &$x_3b_{\ge11}=0$ &  \\
 \hline

\end{tabular}

\end{center}
\caption{The quasi-homogeneous basis of algebraic restrictions of
$2$-forms to the curve $t\mapsto (t^3,t^4,t^5)$.}\label{b345}
\end{table}

The  lowest possible quasi-degree of a $2$-form is $7$. The space
of quasi-homogeneous $2$-forms of degree $7$ is spanned by
$dx_1\wedge dx_2$. This form does not have zero algebraic
restriction since it does not vanish at $0$ \cite{DJZ2}. It
implies that vector space of algebraic restrictions of $2$-forms
of quasi-degree $7$ is spanned by $[dx_1\wedge dx_2]_g$. We have a
similar situation  for the quasi-degrees $8,9,10$. The algebraic
restriction $x_1a_7$ is not zero since there are no
quasi-homogeneous functions vanishing on $g$ of quasi-degree not
greater than $10-3=7$.

The space of quasi-homogeneous $2$-forms of quasi-degree $11$ is
spanned by $x_2dx_1\wedge dx_2$ and $x_1dx_3\wedge dx_1$. But by
Proposition \ref{d-wedge} we have $[df_8\wedge dx_1]_9=0$ which
implies that algebraic restrictions of these $2$-forms are
linearly dependent: $x_1a_8=[x_1dx_3\wedge
dx_1]_g=[-2x_2dx_1\wedge dx_2]_g=-2a_{11}$. We use similar
arguments to show that the space of algebraic restriction of
quasi-degree $12$ is spanned by $a_{12}$.

The space of $2$-forms of quasi-degree $13$ is $3$-dimensional.
But from linearly independent linear relations  satisfied by
algebraic restrictions of elements of the basis presented in the
last column of the row for $\delta=13$ we get that all algebraic
restrictions of quasi-degree $13$ are zero. The same arguments we
use for quasi-degree $14$ and $15$.

To prove that all algebraic restrictions of quasi-degree $16$ are
$0$ we notice that they can have the following forms of
quasi-degree $16$: $x_1\beta_{13}$ or $x_2\beta_{12}$ or
$x_3\beta_{11}$. In the first case the algebraic restriction
$b_{13}=[\beta_{13}]_g$ has quasi-degree $13$, so it is $0$. In
the second case the quasi-degree of $b_{12}=[\beta_{12}]_g$ is
$12$. So the algebraic restriction $b_{12}$ can be presented in
the form $cx_1a_9$, where $c\in \mathbb R$. But then
$x_2b_{12}=x_1(cx_2 a_9)$. The quasi-degree of $cx_2 a_9$ is $13$
and it implies that $cx_2 a_9$ is $0$. We use a similar argument
to prove that $x_3b_{11}$ is $0$. Using the same arguments and
induction by the quasi-degree we show that all algebraic
restrictions of higher quasi-degree are $0$.

Any smooth $2$-form $\omega$ can be decomposed to
$\omega=\sum_{i=7}^{12}\omega_i+\sum_{j=1}^{k}f_j \sigma_j$, where
$k$ is a positive integer, $\omega_i$ is a quasi-homogeneous
$2$-form of quasi-degree $i$ for $i=7,\cdots, 12$ and $f_j$ are
smooth function-germs and $\sigma_j$ are quasi-homogeneous
$2$-forms of quasi-degree greater than $12$ for $j=1,\cdots,k$.
Thus the space of algebraic restrictions of $2$-forms is spanned
by $a_7,\cdots, a_{12}$.

\end{proof}

\begin{prop}\label{bcalr345}
The space of algebraic restrictions of closed differential
$2$-forms to the image of $f$ is the $5$-dimensional vector space
spanned by the following algebraic restrictions:
$$
a_7, \ a_8, \ a_9, \ a_{10}, \ a_{11}.
$$
\end{prop}
\begin{proof}
 It is easy to see that $da_i=0$ for $i<12$ and
 $da_{12}\ne 0$. Then we apply Theorem \ref{th-all-closed}.
 \end{proof}

 \begin{prop}\label{calr345}
 Any algebraic restriction of a symplectic form to $f$ is
 diffeomorphic to one and only one of the following
$a_7, \ a_8, \ -a_8, \ a_9, \ a_{10}, \ -a_{10}, \ a_{11}, \ 0$.
 \end{prop}
\begin{proof} By Theorem \ref{vf-finite} we consider vector fields $X_s$ such that $X_s\circ
f=t^{s+1}df/dt$ for $s=0,...,5$. They have the following form
 $$
 X_0=E=3x_1\frac{\partial}{\partial x_1}+4x_2\frac{\partial}{\partial
 x_2}+5x_3\frac{\partial}{\partial x_3}, \
 X_1=3x_2\frac{\partial}{\partial x_1}+4x_3\frac{\partial}{\partial
 x_2}+5x_1^2\frac{\partial}{\partial x_3},
 $$
 $$
 X_2=3x_3\frac{\partial}{\partial x_1}+4x_1^2\frac{\partial}{\partial
 x_2}+5x_1x_2\frac{\partial}{\partial x_3}, \ X_3=x_1E, \ X_4=x_2E.
 $$

The infinitesimal action of these germs of quasi-homogeneous
liftable vector fields on the basis of the vector space of
algebraic restrictions of closed $2$-forms to $f$ is presented in
Table \ref{ia345}.
\begin{table}
\begin{center}
\begin{tabular}{|c|c|c|c|c|c|}

 \hline

  $\mathcal L_{X_i}a_j$ & $a_7$ & $a_8$ & $a_9$ & $a_{10}$ & $a_{11}$\\ \hline

   $X_0=E$ & $7a_7$ & $8a_8$ & $9a_9$ & $10a_{10}$ & $11a_{11}$\\ \hline

$X_1$ & $-4a_8$ & $-3a_9$ & $-10a_{10}$ & $11a_{11}$ & $0$\\
\hline

$X_2$ & $-3a_9$ & $-5a_{10}$ & $11a_{11}$ & $0$ & $0$\\
\hline

$X_3=x_1 E$ & $10a_{10}$ & $-22a_{11}$ & $0$ & $0$ & $0$\\
\hline

$X_4=x_2 E$ & $11a_{11}$ & $0$ & $0$ & $0$ & $0$\\
\hline

\end{tabular}

\end{center}
\caption{Infinitesimal actions on algebraic restrictions of closed
2-forms to the curve $t\mapsto (t^3,t^4,t^5)$.}\label{ia345}
\end{table}
Using the data of Table \ref{ia345}  we obtain by Theorem
\ref{action} that an algebraic restriction of the form $\sum_{i\ge
s}c_i a_i$ for $c_s\ne 0$ is diffeomorphic to $c_sa_s$. Finally we
reduce $c_sa_s$ to $a_s$ if the quasi-degree $s$ is odd or to
$sgn(c_s)a_s$ if $s$ is even by a diffeomorphism
$\Phi_t(x_1,x_2,x_3)=(t^3x_1,t^4x_2,t^5x_3)$ for
$t=c_s^{\frac{1}{s}}$  or for $t=|c_s|^{\frac{1}{s}}$
respectively.

The algebraic restrictions $a_8$, $-a_8$ are not diffeomorphic.
Any diffeomorphism $\Phi=(\Phi_1,\cdots, \Phi_{2n})$ of $(\mathbb
R^{2n},0)$ preserving $f(t)=(t^3,t^4,t^5,0,\cdots,0)$ has the
following linear part
$$
\begin{array}{ccccccccccccccc}
A^3x_1&+ & A_{12}x_2&+&A_{13}x_3&+&A_{14}x_4&+& \cdots &+&A_{1,2n}x_{2n}\\
& & A^4x_2&+&A_{23}x_3&+&A_{24}x_4&+& \cdots &+&A_{2,2n}x_{2n} \\
& & & & A^5x_3&+&A_{34}x_4&+& \cdots &+&A_{3,2n}x_{2n}\\
& & & & & &A_{44}x_4&+& \cdots &+&A_{4,2n}x_{2n}\\
& & & & & & \vdots& \vdots & \vdots & \vdots & \vdots\\
& & & & & &A_{2n,4}x_4&+& \cdots &+&A_{2n,2n}x_{2n}\\
\end{array}
$$
where $A, A_{i,j}\in \mathbb R$.

Assume that $\Phi^{\ast}(a_8)=-a_8$. It implies that $A^8 dx_3
\wedge dx_1|_0=-dx_3 \wedge dx_1|_0$, which is a contradiction.

One can similarly prove that $a_{10},-a_{10}$ are not
diffeomorphic.
\end{proof}

\begin{proof}[Proof of Theorem \ref{th(3,4,5)}]

Let $\theta_i$ be a $2$-form on $\mathbb R^3$ such that
$a_i=[\theta_i]_g$. Then $rank (\theta_i|_0)\ge 2$ if $n=2$ and
$rank (\theta_i|_0)\ge 0$ if $n>2$ by Proposition \ref{rank}. It
is easy to see that $a_7$, $\pm a_8$, $a_9$ are realizable by the
following  symplectic forms
$$
dx_1\wedge dx_2+dx_3\wedge dx_4+\cdots, \ \pm dx_3\wedge
dx_1+dx_2\wedge dx_4+\cdots, \ dx_2\wedge dx_3+dx_1\wedge
dx_4+\cdots
$$
respectively. The algebraic restrictions $\pm a_{10}$, $a_{11}$
and $a_{\infty}=0$ are realizable by the following forms
$$\pm x_1dx_1\wedge dx_2+dx_1\wedge dx_4
+dx_2\wedge dx_5 +dx_3 \wedge dx_6+\cdots,$$ $$ x_2dx_1\wedge
dx_2+dx_1\wedge dx_4 +dx_2\wedge dx_5 +dx_3 \wedge dx_6+\cdots,$$
$$dx_1\wedge dx_4 +dx_2\wedge dx_5 +dx_3 \wedge dx_6+\cdots$$
respectively. By a simple coordinate change we map the above forms
to the Darboux normal form and we obtain the normal forms of the
curve.

By Propositions \ref{sm}, \ref{bcalr345}, \ref{calr345} and using
the data in Table \ref{ia345} we obtain the symplectic
multiplicities of curves in Table \ref{ss345}.  The indexes of
isotropness for these curves are calculated by Propositions
\ref{ii} and \ref{calr345}. The Lagrangian tangency orders for the
curves in rows $1-3$ are obtained using the fact that any
Lagrangian submanifold can be represented in the form (\ref{L}).
By Propositions \ref{lto} and \ref{calr345} we obtain this
invariant for other curves in Table \ref{ss345}.

\end{proof}

\section{Symplectic singularities of curves with the semigroup $(3,5,7)$}\label{357}
In this section we present the symplectic classification of curves
with the semigroup $(3,5,7)$.
\begin{thm} Let $(\mathbb
R^{2n},\omega_0=\sum_{i=1}^n dp_i\wedge dq_i$) be the symplectic
space  with the canonical coordinates $(p_1,q_1,\cdots, p_n,q_n)$.

Then the germ of a curve $f:(\mathbb R,0)\rightarrow (\mathbb
R^{2n},0)$ with the semigroup $(3,5,7)$ is symplectically
equivalent to one and only one of the curves presented in the
second column of the Table \ref{ss357} (on page 16) for $n>2$ and
$f$ is symplectically equivalent to one and only one of the curves
presented in the second column and rows 1-3 and 5 for $n=2$. The
parameter $c$ is a modulus.

The symplectic multiplicity, the index of isotropness and the
Lagrangian tangency order are presented in the third, fourth and
fifth columns of Table \ref{ss357}.
\begin{table}[h]
\begin{center}
\begin{tabular}{|c|l|c|c|c|}

 \hline

   & normal form of $f$ & $\mu_{sympl}(f)$& $\iota(f)$ & $Lt(f)$\\ \hline

   1 &$t\mapsto(t^3,\pm t^5,t^7,0,\cdots,0)$ & $0$ & $0$& $5$\\ \hline

 2 &$t\mapsto(t^3,\pm t^7,t^5,ct^6,\cdots,0)$&$2$  &$0$  &$7$\\ \hline

 3 & $t\mapsto(t^3,t^8,t^5,ct^7,\cdots,0)$, $c\ne 0$ & $3$ &$0$  &$7$\\ \hline

4 &$t\mapsto(t^3,t^8,t^5,0,t^7,0,\cdots,0)$ &$3$  &$1$  &$8$\\
\hline

5 &$t\mapsto(t^3,ct^{10},\pm t^5,t^7,0,\cdots,0)$ &$4$  &$0$  &$7$\\
\hline

6 & $t\mapsto(t^3,ct^{11},t^5,t^8,t^7,0,\cdots,0)$ & $5$ &$1$  &$10$\\
\hline

7 & $t\mapsto(t^3,\pm t^{11},t^5,0,t^7,0,\cdots,0)$ & $5$ &$2$  &$11$\\
\hline

8 & $t\mapsto(t^3,\pm t^{13}/2,t^5,0,t^7,0,\cdots,0)$ & $6$ &$2$  &$13$\\
\hline

9 & $t\mapsto(t^3,0,t^5,0,t^7,0,\cdots,0)$ & $7$ &$\infty$  &$\infty$\\
\hline
\end{tabular}
\end{center}
\caption{Symplectic classification of curves with the semigroup
$(3,5,7)$.}\label{ss357}
\end{table}

\end{thm}

The germ of a curve $f:\mathbb R\ni t\mapsto f(t)\in \mathbb
R^{2n}$ with the semigroup $(3,5,7)$ is $\mathcal A$-equivalent to
$t\mapsto (t^3,t^5,t^7,0,\cdots,0)$. We use the same method as in
the previous section to obtain symplectic classification of these
curves. We only present the main steps with all calculation
results in tables.
\medskip
\begin{prop}\label{basis-all-2}
The space of algebraic restrictions of differential $2$-forms to
$g$ is the $8$-dimensional vector space spanned by the following
algebraic restrictions:
$$
a_8=[dx_1\wedge dx_2]_g,  a_{10}=[dx_3\wedge dx_1]_g,
a_{11}=[x_1dx_1\wedge dx_2]_g,  a_{12}=[dx_2\wedge dx_3]_g,
$$
$$
a_{13}=[x_2 dx_1\wedge dx_2]_g,  a_{14}=[x_1^2dx_1\wedge dx_2]_g,
a_{15}=[x_3dx_1\wedge dx_2]_g,  a_{16}=[x_1x_2dx_1\wedge dx_2]_g,
$$
where $\delta$ is the  quasi-degree of $a_{\delta}$.
\end{prop}

\begin{proof}[The sketch of the proof]
We use the same method as in the previous section. The sketch of
the proof is presented in Tables \ref{v357} and \ref{ab357}.
\end{proof}
\begin{table}
\begin{center}
\begin{tabular}{|c|c|c|}

 \hline

  quasi-degree $\delta$ & $h_{\delta}$& differential $dh_{\delta}$\\ \hline

   $10$ & $x_1x_3-x_2^2$ & $x_1dx_3+x_3dx_1-2x_2dx_2$\\ \hline

  $12$ & $x_2x_3-x_1^4$ & $x_2dx_3+x_3dx_2-4x_1^3dx_1$\\ \hline

$14$ & $x_1^3x_2-x_3^2$ & $3x_1^2x_2dx_1+x_1^3dx_2-2x_3dx_3$\\
\hline
\end{tabular}
\end{center}
\caption{Quasi-homogeneous function-germs of quasi-degree $10, \
12, \ 14$ vanishing on the curve $t\mapsto
(t^3,t^5,t^7)$.}\label{v357}
\end{table}
\begin{prop}
The space of algebraic restrictions of closed differential
$2$-forms to the image of $f$ is the $7$-dimensional vector space
spanned by the following algebraic restrictions:
$$
a_8, \ a_{10}, \ a_{11}, \ a_{12}, \ a_{13}, \ a_{14}, \ a_{16}.
$$
\end{prop}
\begin{proof}
By Proposition \ref{basis-all-2}
 it is easy to see that $da_{i}=0$ for $i\ne 15$ and
 $da_{15}\ne 0$. By Theorem \ref{th-all-closed} we get the result.
\end{proof}

 \begin{prop}
 Any algebraic restriction of a symplectic form to $f$ is
 diffeomorphic to one of the following
$\pm a_8, \ \pm a_{10}+ca_{11}, \ a_{11}+ca_{12}, \ a_{11}, \ \pm
a_{12}+ca_{13}, \ a_{13}+ca_{14}, \ \pm a_{14}, \ \pm a_{16},\ 0$,
where the parameter $c\in \mathbb R$ is a modulus.
 \end{prop}
\begin{proof}[Sketch of the proof] The vector fields $X_s$ (see Theorem
\ref{vf-finite}) have the following form:
$$
 X_0=E=3x_1\frac{\partial}{\partial x_1}+5x_2\frac{\partial}{\partial
 x_2}+7x_3\frac{\partial}{\partial x_3}, \
 X_2=3x_2\frac{\partial}{\partial x_1}+5x_3\frac{\partial}{\partial
 x_2}+7x_1^3\frac{\partial}{\partial x_3},
 $$
 $$
 X_3=x_1E, \
 X_4=3x_3\frac{\partial}{\partial x_1}+5x_1^2\frac{\partial}{\partial
 x_2}+7x_1^2x_2\frac{\partial}{\partial x_3},
 $$
 $$
 X_5=x_2E, \ X_6=x_1^2E, \
 X_7=x_3E, \ X_8=x_1x_2E.
 $$
Their actions on the space of algebraic restrictions of closed
$2$-forms are presented in Table \ref{ia357}. From these data we
obtain the classification of algebraic restrictions as in the
previous section.

From Table \ref{ia357} and Theorem \ref{vf-finite} we also see
that the tangent space to the orbit of $\pm a_{10} + ca_{11}$ at
$\pm a_{10} + ca_{11}$ is spanned by $\pm 10a_{10}+11c a_{11}$,
$a_{12}$, $a_{13}$, $a_{14}$, $a_{16}$. $a_{11}$ does not belong
to it. Therefore parameter $c$ is the modulus in the  normal form
$\pm a_{10} + ca_{11}$.

In the same way we prove that $c$ is the modulus in  the other
normal forms.
 \end{proof}

\begin{table}
\begin{center}
\begin{tabular}{|c|c|c|c|c|}

 \hline

  $\delta$&  basis & forms & relations & proof\\ \hline

   $8$ & $a_8$&$\alpha_8=dx_1\wedge dx_2$&$a_8:=[\alpha_8]_g$ &\\ \hline

  $10$ & $a_{10}$ & $\alpha_{10}=dx_3\wedge dx_1$&$a_{10}:=[\alpha_{10}]_g$ &\\ \hline

$11$ &  $a_{11}$ & $x_1\alpha_8$& $a_{11}:=x_1a_8$&\\
\hline

$12$ &  $a_{12}$& $\alpha_{12}=dx_2\wedge dx_3$ &$a_{12}:=[\alpha_{12}]_g$  &\\
\hline

$13$ & $a_{13}$ & $x_2\alpha_8, $ & $a_{13}:=x_2a_8$ & $[dh_{10}\wedge dx_1]_g=0$\\
& &$x_1\alpha_{10}$ & $x_1a_{10}=-2a_{13}$ & \\
\hline

$14$ &  $a_{14}$ & $x_1^2\alpha_8$& $a_{14}:=x_1^2a_8$&\\
\hline

 $15$ & $a_{15}$ &$x_3\alpha_8,$ &
$a_{15}:=x_3a_8$   & $[dh_{12}\wedge dx_1]_g=0$\\
& &$x_2\alpha_{10}, $ &$x_2a_{10}=a_{15}$ & $[dh_{10}\wedge
dx_2]_g=0$ \\
& &$x_1 \alpha_{12}$ &$x_1a_{12}=a_{15}$ &\\
\hline

$16$ & $a_{16}$ & $x_1x_2\alpha_8, $ & $a_{16}:=x_1x_2a_8$ & $x_1[dh_{10}\wedge dx_1]_g=0$\\
& &$x_1^2\alpha_{10}$ & $x_1^2a_{10}=-2a_{16}$ & \\
\hline

$17$ & $0$ &$x_1^3\alpha_8, $  &$x_1^3a_8=0$
  & $[dh_{14}\wedge dx_1]_g=0$\\
& &$x_3\alpha_{10},$ & $x_3a_{10}=0$& $[dh_{12}\wedge dx_2]_g=0$ \\
& &$x_2\alpha_{12}$ & $x_2a_{12}=0$& $[dh_{10}\wedge dx_3]_g=0$ \\
\hline

$18$ & $0$ &$x_1x_3\alpha_8,$ & $x_1x_3a_8=0$& $x_1[dh_{12}\wedge dx_1]_g=0$ \\
& &$x_1x_2\alpha_{10}$ & $x_1x_2a_{10}=0$& $x_1[dh_{10}\wedge
dx_2]_g=0$ \\
& &$x_1^2\alpha_{12}$ & $x_1^2a_{12}=0$& $x_2[dh_{10}\wedge dx_1]_g=0$ \\
\hline

$19$ & $0$ &$x_1^2x_2\alpha_8, $  &$x_1^2x_2a_8=0$ & $x_1^2[dh_{10}\wedge dx_1]_g=0$ \\
& & $x_1^3\alpha_{10}, $ & $x_2^2a_{10}=0$
  & $[dh_{12}\wedge dx_3]_g=0$\\
& &$x_3\alpha_{12},$ & $x_3a_{12}=0$& $[dh_{14}\wedge dx_2]_g=0$ \\
\hline

$\ge 20$ & $0$ &$x_1\beta_{\ge17}, $ &$
b_{\ge17}:=[\beta_{\ge17}]_g$
  & $x_2b_{\ge15}=x_1b'_{\ge 17}$\\
& & &$x_1b_{17}=0$& $x_3b_{\ge13}=x_1b''_{\ge 17}$\\
& &$x_2\beta_{\ge15},$ &$b_{\ge15}=:[\beta_{\ge15}]_g$ & $\delta(b_{\ge 17})\ge 17$ \\
& & &$x_2b_{\ge15}=0$ & $ b_{\ge 17}=0$ \\
& &$x_3\beta_{\ge13}$ & $b_{\ge13}:=[\beta_{\ge13}]_g$&  \\
& & &$x_3b_{\ge13}=0$ &  \\

\hline
\end{tabular}
\end{center}
\caption{The quasi-homogeneous basis of algebraic restrictions of
$2$-forms  to the curve $t\mapsto (t^3,t^5,t^7)$.}\label{ab357}
\end{table}
\begin{table}
\begin{center}
\begin{tabular}{|c|c|c|c|c|c|c|c|}

 \hline

  $\mathcal L_{X_i}a_j$ & $a_8$ & $a_{10}$ & $a_{11}$ & $a_{12}$ & $a_{13}$ & $a_{14}$ & $a_{16}$ \\ \hline

   $X_0=E$ & $8a_8$ & $10a_{10}$ & $11a_{11}$ & $12a_{12}$ & $13a_{13}$ & $14a_{14}$ & $16a_{16}$ \\ \hline

$X_2$ & $-5a_{10}$ & $-3a_{12}$ & $13a_{13}$ & $-21a_{14}$ & $0$ & $16a_{16}$ & $0$ \\
\hline

$X_3=x_1E$ & $11a_{11}$ & $-26a_{13}$ & $14a_{14}$ & $0$ & $16a_{16}$ & $0$ & $0$ \\
\hline

$X_4$ & $-3a_{12}$ & $-7a_{14}$ & $0$ & $6a_{16}$ & $0$ & $0$ & $0$ \\
\hline

$X_5=x_2 E$ & $13a_{13}$ & $0$ & $16a_{16}$ & $0$ & $0$ & $0$ & $0$ \\
\hline

$X_6=x_1^2 E$ & $11a_{14}$ & $-32a_{16}$ & $0$ & $0$ & $0$ & $0$ & $0$ \\
\hline

$X_7=x_3 E$ & $0$ & $0$ & $0$ & $0$ & $0$ & $0$ & $0$ \\
\hline

$X_8=x_1x_2 E$ & $4a_{16}$ & $0$ & $0$ & $0$ & $0$ & $0$ & $0$ \\
\hline
\end{tabular}
\end{center}
\caption{Infinitesimal actions on algebraic restrictions of closed
2-forms to the curve $t\mapsto (t^3,t^5,t^7)$.}\label{ia357}
\end{table}

\section{Symplectic singularities of curves with the semigroup $(3,7,8)$}\label{378}
In this section we present the symplectic classification of curves
with the semigroup $(3,7,8)$.
\begin{thm} Let $(\mathbb
R^{2n},\omega_0=\sum_{i=1}^n dp_i\wedge dq_i$) be the symplectic
space  with the canonical coordinates $(p_1,q_1,\cdots, p_n,q_n)$.

Then the germ of a curve $f:(\mathbb R,0)\rightarrow (\mathbb
R^{2n},0)$ with the semigroup $(3,7,8)$ is symplectically
equivalent to one and only one of the curves presented in the
second column of the Table \ref{ss378} (on page 19) for $n>2$ and
$f$ is symplectically equivalent to one and only one of the curves
presented in the second column and rows 1-3, 5 and 7 for $n=2$.
The parameters $c$, $c_1$, $c_2$ are moduli.

The symplectic multiplicity, the index of isotropness and the
Lagrangian tangency order are presented in the third, fourth and
fifth columns of Table \ref{ss378}.
\begin{table}[h]
\begin{center}
\begin{tabular}{|c|l|c|c|c|}

 \hline

   & normal form of $f$ & $\mu_{sympl}(f)$& $\iota(f)$ & $Lt(f)$\\ \hline

   1 &$t\mapsto(t^3,\pm t^7,t^8,ct^3,0,\cdots,0)$ & $1$ & $0$& $7$\\ \hline

 2 &$t\mapsto(t^3,t^8,t^7,ct^6,0,\cdots,0)$&$2$  &$0$  &$8$\\ \hline

 3 & $t\mapsto(t^3,t^{10}+c_1t^{11},t^7,c_2t^8,0,\cdots,0)$, $c_2\ne 0$ & $4$ &$0$  &$8$\\ \hline

4 &$t\mapsto(t^3,t^{10},t^8,ct^6,t^7,0,\cdots,0)$ &$4$  &$1$  &$10$\\
\hline

5 &$t\mapsto(t^3,\pm t^{11}+c_2t^{13},t^7,c_1t^8,0,\cdots,0)$, $c_1\ne 0$ &$5$  &$0$  &$8$\\
\hline

6 & $t\mapsto(t^3,\pm t^{11},t^7,ct^9,t^8,0,\cdots,0)$
 & $5$ &$1$  &$11$\\
\hline

7 & $t\mapsto(t^3,c_1 t^{13}+c_2
t^{14},t^7,t^8,0,\cdots,0)$ & $6$ &$0$  &$8$\\
\hline

8 & $t\mapsto(t^3,\pm
t^{13},t^7,c_1t^{10},t^8,c_2t^{10},0,\cdots,0)$, $c_2\ne 0$& $7$ &$1$  &$11$  \\
\hline

9 & $t\mapsto(t^3,\pm t^{13},t^7,ct^{10},t^8,0,\cdots,0)$& $7$
&$2$ & $13$ \\ \hline

10 &
$t\mapsto(t^3,t^{14},t^7,c_1t^{11},t^8,c_2t^{11},0,\cdots,0)$,
$c_1\ne 0$ & $8$ &$1$ &$11$\\ \hline
11 &$t\mapsto(t^3,t^{14},t^7,0,t^8,ct^{11},0,\cdots,0)$ &$8$ &$2$ & $14$\\
\hline

12 & $t\mapsto(t^3,c_1t^{17},t^7,\pm
t^{11},t^8,c_2t^{11},0,\cdots,0)$ & $9$ &$1$  &$11$\\
\hline

13 & $t\mapsto(t^3,t^{16},t^7,c t^{13},t^8,0,\cdots,0)$
 & $9$ &$3$  &$16$\\
\hline

14 & $t\mapsto(t^3,\pm t^{17},t^7,0,t^8,0,\cdots,0)$ & $9$ &$3$  &$17$\\
\hline

15 & $t\mapsto(t^3,0,t^7,0,t^8,0,\cdots,0)$
 & $10$ &$\infty$  &$\infty$\\
\hline
\end{tabular}

\end{center}
\caption{Symplectic classification of curves with the semigroup
$(3,7,8)$.}\label{ss378}
\end{table}

\end{thm}
Let $f:\mathbb R\ni t\mapsto f(t)\in \mathbb R^{2n}$ be the germ
of a smooth or $\mathbb R$-analytic curve $\mathcal A$-equivalent
to $t\mapsto (t^3,t^7,t^8,0,\cdots,0)$. First we calculate the
space of algebraic restrictions of $2$-forms to the image of $f$
in $\mathbb R^{2n}$.

\begin{prop}\label{basis-all-378}
The space of algebraic restrictions of differential $2$-forms to
$g$ is the $12$-dimensional vector space spanned by the following
algebraic restrictions:
$$
a_{10}=[dx_1\wedge dx_2]_g,  a_{11}=[dx_3\wedge dx_1]_g,
a_{13}=[x_1dx_1\wedge dx_2]_g,  a_{14}=[x_1dx_3\wedge dx_1]_g,
$$
$$
a_{15}=[dx_2\wedge dx_3]_g, a_{16}=[x_1^2dx_1\wedge dx_2]_g,
a_{17}=[x_2 dx_1\wedge dx_2]_g, a_{18}^+=[x_1dx_2\wedge dx_3]_g,
$$
$$
a_{18}^-=[x_2dx_3\wedge dx_1]_g, a_{19}=[x_3dx_3\wedge dx_1]_g,
a_{20}=[x_1x_2dx_1\wedge dx_2]_g, a_{21}=[x_1x_3dx_1\wedge
dx_2]_g,
$$
where $\delta$ is quasi-degree of $a_{\delta}$.
\end{prop}

\begin{proof}[The sketch of  the proof]
We use the same method as in the previous sections. The sketch of
the proof is presented in Tables \ref{v378} and \ref{ab378}.
\end{proof}
\begin{table}
\begin{center}
\begin{tabular}{|c|c|c|}

 \hline

  quasi-degree $\delta$ & $h_{\delta}$& differential $dh_{\delta}$\\ \hline

   $14$ & $x_1^2x_3-x_2^2$ & $2x_1x_3dx_1+x_1^2dx_3-2x_2dx_2$\\ \hline

  $15$ & $x_2x_3-x_1^5$ & $x_2dx_3+x_3dx_2-5x_1^4dx_1$\\ \hline

$16$ & $x_1^3x_2-x_3^2$ & $3x_1^2x_2dx_1+x_1^3dx_2-2x_3dx_3$\\
\hline
\end{tabular}
\end{center}
\caption{Quasi-homogeneous function-germs of quasi-degree $14, \
15, \ 16$ vanishing on the curve $t\mapsto
(t^3,t^7,t^8)$.}\label{v378}
\end{table}
\begin{table}
\begin{center}
\begin{tabular}{|c|c|c|c|c|}

 \hline

  $\delta$&  basis & forms & relations & proof\\ \hline

   $10$ & $a_{10}$&$\alpha_{10}=dx_1\wedge dx_2$&$a_{10}:=[\alpha_{10}]_g$ &\\ \hline

  $11$ & $a_{11}$ & $\alpha_{11}=dx_3\wedge dx_1$&$a_{11}:=[\alpha_{11}]_g$ &\\ \hline

$13$ &  $a_{13}$ & $x_1\alpha_{10}$& $a_{13}:=x_1a_{10}$&\\
\hline

$14$ &  $a_{14}$& $\alpha_{14}=x_1\alpha_{11}$ &$a_{14}:=x_1a_{11}$  &\\
\hline

$15$ & $a_{15}$ & $\alpha_{15}=dx_2\wedge dx_3 $ &
$a_{15}:=[\alpha_{15}]_g$ &\\ \hline

$16$ &  $a_{16}$ & $x_1^2\alpha_{10}$& $a_{16}:=x_1^2a_{10}$&\\
\hline

$17$ &  $a_{17}$ & $x_2\alpha_{10},$&
$a_{17}:=x_2a_{10}$&$[dh_{14}\wedge
dx_1]_g=0$\\
& &$x_1^2\alpha_{11} $ &$x_1^2a_{11}=-2a_{17}$ &  \\
\hline

 $18$ & $a^+_{18},$ &$x_1\alpha_{15},$ &
$a^+_{18}:=x_1a_{15}$   & $[dh_{15}\wedge dx_1]_g=0$\\
&$a^-_{18}$ &$x_2\alpha_{11}, $ &$a^-_{18}:=x_2a_{11}$ &
 \\
& &$x_3 \alpha_{10}$ &$x_3a_{10}=a^-_{18}$ &\\
\hline

$19$ & $a_{19}$ & $x_3\alpha_{11}, $ & $a_{19}:=x_3a_{11}$ & $[dh_{16}\wedge dx_1]_g=0$\\
& &$x_1^3\alpha_{10}$ & $x_1^3a_{10}=-2a_{19}$ & \\
\hline

$20$ & $a_{20}$ &$x_1x_2\alpha_{10}, $  &$a_{20}:=x_1x_2a_{10}$
  & $x_1[dh_{14}\wedge dx_1]_g=0$\\
& &$x_1^3\alpha_{11}$ & $x_1^3a_{10}=-2a_{20}$& \\
\hline

$21$ & $a_{21}$ &$x_1x_3\alpha_{10}, $  &$a_{21}:=x_1x_3a_{10}$
  & $[dh_{14}\wedge dx_2]_g=0$\\
& &$x_1x_2\alpha_{11},$ & $x_1x_2a_{11}=a_{21}$& $x_1[dh_{15}\wedge dx_1]_g=0$ \\
& &$x_1^2\alpha_{15}$ & $x_1^2a_{15}=2a_{21}$&  \\
\hline

$22$ & $0$ &$x_1^4\alpha_{10}, $  &$x_1^4a_{10}=0$ & $x_1[dh_{16}\wedge dx_1]_g=0$ \\
& & $x_1x_3\alpha_{11}, $ & $x_1x_3a_{11}=0$
  & $[dh_{15}\wedge dx_2]_g=0$\\
& &$x_2\alpha_{15},$ & $x_2a_{15}=0$& $[dh_{14}\wedge dx_3]_g=0$ \\
\hline

$23$ & $0$ &$x_1^2x_2\alpha_{10}, $  &$x_1^2x_2a_{10}=0$ & $x_1^2[dh_{14}\wedge dx_1]_g=0$ \\
& & $x_1^4\alpha_{11}, $ & $x_1^4a_{11}=0$
  & $[dh_{16}\wedge dx_2]_g=0$\\
& &$x_3\alpha_{15},$ & $x_3a_{15}=0$& $[dh_{15}\wedge dx_3]_g=0$ \\
\hline

$24$ & $0$ &$x_2^2\alpha_{10}, $  &$x_2^2a_{10}=0$ & $x_1^2[dh_{15}\wedge dx_1]_g=0$ \\
& & $x_1^2x_2\alpha_{11}, $ & $x_1^2x_2a_{11}=0$
  & $[dh_{16}\wedge dx_3]_g=0$\\
& &$x_3\alpha_{15},$ & $x_3a_{15}=0$& $x_1[dh_{14}\wedge dx_2]_g=0$ \\
\hline

$\ge 25$ & $0$ &$x_1\beta_{\ge22}, $ &$
b_{\ge22}:=[\beta_{\ge22}]_g$
  & $x_2b_{\ge 18}=x_1b'_{\ge 22}$\\
& & &$x_1b_{\ge22}=0$& $x_3b_{\ge 17}=x_1b''_{\ge 22}$\\
& &$x_2\beta_{\ge18},$ &$b_{\ge18}=:[\beta_{\ge18}]_g$ & $ \delta(b_{\ge 22})\ge 22$ \\
& & &$x_2b_{\ge18}=0$ & $ b_{\ge 22}=0 $ \\
& &$x_3\beta_{\ge17}$ & $b_{\ge17}:=[\beta_{\ge17}]_g$&  \\
& & &$x_3b_{\ge17}=0$ &  \\
\hline
\end{tabular}
\end{center}
\caption{The quasi-homogeneous basis of algebraic restrictions of
$2$-forms to the curve $t\mapsto (t^3,t^7,t^8)$.}\label{ab378}
\end{table}
\begin{prop}
The space of algebraic restrictions of closed differential
$2$-forms to the image of $g$ is the $10$-dimensional vector space
spanned by the following algebraic restrictions:
$$
a_{10}, \ a_{11}, \ a_{13}, \ a_{14}, \ a_{15}, \ a_{16}, \
a_{17}, \ a_{18}=a_{18}^+ - a_{18}^-, \ a_{19}, \ a_{20}
$$
\end{prop}
\begin{proof}

 It is easy to see that
 $da_i=0$ for $i\ne 18,21$,
 $da_{18}^+=da_{18}^-\ne 0$ and $da_{21}\ne 0$.
 Then the algebraic restriction
 $a_{18}^{+}-a_{18}^{-}$ is closed and $da^-_{18}$, $da_{21}$ are linearly independent.
 Thus  Theorem \ref{th-all-closed}
 implies the result.
 \end{proof}

 \begin{prop}
 Any algebraic restriction of a symplectic form to $f$ is
 diffeomorphic to one of the following
$\pm a_{10}+ca_{11}, \ a_{11}+ca_{13}, \
a_{13}+c_1a_{14}+c_2a_{15}, \ \pm a_{14}+c_1a_{15}+c_2a_{16}, \
a_{15}+c_1a_{16}+c_2a_{17},\ \pm a_{16}+c_1a_{17}+c_2a_{18},\
a_{17}+c_1a_{18}+c_2a_{19}, \ \pm a_{18}+c_1a_{19}+c_2a_{20}, \
a_{19}+ca_{20},\ \pm a_{20},\ 0$, where $c, c_1, c_2 \in \mathbb
R$.

The parameters $c$, $c_1$, $c_2$ are moduli.
 \end{prop}

\begin{proof} [The sketch of the proof] The vector fields $X_s$ (see Theorem
\ref{vf-finite}) have the following form:

 $$
 X_0=E=3x_1\frac{\partial}{\partial x_1}+7x_2\frac{\partial}{\partial
 x_2}+8x_3\frac{\partial}{\partial x_3}, \ X_3=x_1E,
 $$
 $$
 X_4=3x_2\frac{\partial}{\partial x_1}+7x_1x_3\frac{\partial}{\partial
 x_2}+8x_1^4\frac{\partial}{\partial x_3}, \ X_5=3x_3\frac{\partial}{\partial x_1}+7x_1^4\frac{\partial}{\partial
 x_2}+8x_1^2x_2\frac{\partial}{\partial x_3},
 $$
 $$
 X_6=x_1^2E, \ X_7=x_2E, \
 X_8=x_3E, \ X_9=x_1^3E,
 X_{10}=x_1x_2E.
 $$
Their actions on the space of algebraic restrictions of closed
$2$-forms are presented in Table \ref{ia378}. From these data we
obtain the classification of algebraic restrictions as in the
previous section.

Now we prove that parameters $c$, $c_1$, $c_2$ are moduli in the
normal forms. The proofs are very similar in all cases. As an
example we consider the normal form $a_{13}+c_1a_{14}+c_2a_{15}$ -
the first normal form with two parameters. From Table \ref{ia378}
and Theorem \ref{vf-finite} we see that the tangent space to the
orbit of $a_{13}+c_1a_{14}+c_2a_{15}$ at
$a_{13}+c_1a_{14}+c_2a_{15}$ is spanned by linearly independent
algebraic restrictions $a_{13}+14c_1a_{14}+15c_2a_{15}$, $a_{16}$,
$a_{17}$, $a_{18}$, $a_{19}$, $a_{20}$. Hence algebraic
restrictions $a_{14}$ and $a_{15}$ do not belong to it. Therefore
parameters $c_1$ and $c_2$ are independent moduli in the normal
form $a_{13}+c_1a_{14}+c_2a_{15}$.
\end{proof}

\begin{table}
\begin{center}
\begin{tabular}{|c|c|c|c|c|c|c|c|c|}

 \hline

  $\mathcal L_{X_i}a_j$  & $a_{10}$ & $a_{11}$ & $a_{13}$ & $a_{14}$ & $a_{15}$ &$a_{16}$ & $a_{17}$& $a_{\delta\ge 18}$\\ \hline

   $E$ & $10a_{10}$ & $11a_{11}$ & $13a_{13}$ & $14a_{14}$ & $15a_{15}$ & $16a_{16}$& $17a_{17}$& $\delta a_{\delta}$\\ \hline

$X_3$ & $13a_{13}$ & $14a_{14}$ & $16a_{16}$ & $17a_{17}$ & $15a_{18}$ & $-38a_{19}$ & $20a_{20}$ & $0$ \\
\hline

$X_4$ & $-7a_{14}$ & $-3a_{15}$ & $17a_{17}$ & $-3a_{18}$ & $9a_{19}$ & $20a_{20}$ & $0$ & $0$\\
\hline

$X_5$ & $-3a_{15}$ & $-8a_{16}$ & $-3a_{18}$ & $19a_{19}$ & $12a_{20}$ & $0$ & $0$ & $0$\\
\hline

$X_6$ & $16a_{16}$ & $17a_{17}$ & $-38 a_{19}$ & $-38a_{20}$ & $0$ & $0$ & $0$ & $0$\\
\hline

$X_7$ & $17a_{17}$ & $-3a_{18}$ & $20a_{20}$ & $0$ & $0$ & $0$ & $0$ & $0$\\
\hline

$X_8$ & $-3a_{18}$ & $19a_{19}$ & $0$ & $0$ & $0$ & $0$ & $0$ & $0$\\
\hline

$X_9$ & $-38a_{19}$ & $-40a_{20}$ & $0$ & $0$ & $0$ & $0$ & $0$ & $0$\\
\hline

$X_{10}$ & $20a_{20}$ & $0$ & $0$ & $0$ & $0$ & $0$ & $0$ & $0$\\
\hline
\end{tabular}
\end{center}
\caption{Infinitesimal actions on algebraic restrictions of closed
2-forms to the curve $t\mapsto (t^3,t^7,t^8)$.}\label{ia378}
\end{table}

\bibliographystyle{amsalpha}

\end{document}